\documentclass[11pt,reqno]{amsart}

\usepackage{amsmath, amsthm, amssymb,amsaddr}
\usepackage{dcolumn}
\usepackage{bm}
\usepackage{amsmath} 
\usepackage{mathrsfs}
\usepackage{mathtools}
\usepackage{algpseudocode}
\usepackage{algorithm}
\usepackage[pdftex]{graphicx}
\usepackage{subfigure}
\usepackage{epstopdf}
\usepackage{tikz}
\usepackage{xcolor}
\usepackage[export]{adjustbox}
\usepackage{caption}
\usepackage{hyperref}
\usepackage{stmaryrd}
\usepackage[utf8]{inputenc} 
\usepackage[T1]{fontenc}

\theoremstyle{remark}
\newtheorem{remark}{Remark}

\newtheorem{definition}{Definition}

\newtheoremstyle{italic_lemma} 
{} 
{} 
{\itshape} 
{} 
{\itshape} 
{.} 
{ } 
{} 

\theoremstyle{italic_lemma}
\newtheorem{lemma}{Lemma}

\DeclareMathOperator{\sinc}{sinc}

\newcommand{\eps}{\varepsilon}

\newcommand{\be}[1]{\begin{equation}\label{#1}}
\newcommand{\ee}{\end{equation}}

\newcommand{\no}[1]{#1}

\newcommand{\bt}{\mathfrak{t}}
\newcommand{\bx}{\boldsymbol{x}}

\newcommand{\bu}{\boldsymbol{u}}
\newcommand{\bro}{\boldsymbol{\rho}}
\newcommand{\bX}{\boldsymbol{X}}

\newcommand{\bn}{\boldsymbol{n}}

\newcommand{\fd}{\mathfrak{d}}
\newcommand{\btx}{\tilde{\bx}}

\renewcommand{\no}[1]{} 

\newcounter{subsubsubsection}[subsubsection]

\title[]
{A Regularized Finite-Difference Approximation of Surface-restricted Emission and Reception Process in Acoustics with Application to Inverse Problems}
\author{Ashkan Javaherian}
\email{ajavaherian62@gmail.com; ashkan.javaherian@ut.ac.ir}
\address{Department of Bio-Electric, School of Electrical and Computer Engineering, \\
  University College of Engineering, University of Tehran, Tehran, Iran.}
\date{September 2026}

\begin{document}
\maketitle

\section*{Abstract} 
\noindent
The acoustic wave equation governs wave propagation induced by either volumetric radiation sources, or by surface sources of monopole or dipole type. For surface sources, boundary value problems yield wavefield representations via the Kirchhoff–Helmholtz or Rayleigh–Sommerfeld integrals. This study begins by examining the equivalence between the analytic expressions of the associated monopole and dipole integral formulations and their regularized approximations. Leveraging these regularized formulations, we introduce \emph{reception operators} that map free space pressure wavefields—obtained by solving the wave equation—onto measured fields restricted to the boundary. Building on this trace mapping, we derive the adjoint of the forward operator. We show that, under the common practical assumption of Dirichlet-type boundary data, the adjoint operator coincides—up to a constant factor—with the time-reversed form of the dipole integral formula, evaluated on the receiver surfaces. This study aims to advance the numerical approximation of forward problems and the solution of inverse problems in acoustics, with a particular focus on applications that require accurate amplitude modeling, including attenuation reconstruction and photoacoustic tomography. \footnote{An extension of the findings of this study to the inverse problem of photoacoustic tomography has been reported in \cite{Javaherian10}.}

\section{Introduction}
\noindent
The acoustic wave equation is one of the most important partial differential equations (PDEs) in mechanics, with broad applications across a wide range of fields~\cite{Mast,Tabei,Treeby,Bilbao,Bencomo,Borcea,Kaltenbacher}.

The solution to the wave equation is typically non-unique. To obtain a unique solution, Cauchy initial conditions, defined in terms of the wavefield and its time derivative at the initial time, are enforced. These initial conditions establish a causal relationship between the solution wavefield and the source, ensuring that the solution wavefield vanishes prior to the onset of radiation from the source. Thus, solving the wave equation can be formulated as uniquely representing the propagated wavefield in terms of the radiation source (forcing term) by applying the causality conditions \cite{Devaney}.


 

Let $d \in \{2,3\}$ denote the number of spatial dimensions of the medium. The radiation source, $s$, is defined over a finite $d$-dimensional space (a volumetric region for $d=3$) and time, and is assumed to possess finite energy; that is, it is square-integrable over the full space-time domain. Furthermore, the wavefield may also be represented in terms of a surface-restricted source~\cite{Wu1}, where the surface either bounds the region supporting the radiation source or lies on an infinite plane. The wavefield induced by a source confined to a surface can be described using the \textit{Kirchhoff--Helmholtz} or \textit{Rayleigh--Sommerfeld} integral formulas, where the integration is carried out over this surface, referred to here as the \textit{acoustic aperture}. The former integral is used for apertures with arbitrary geometries, such as spherical or cylindrical surfaces, whereas the latter represents a special case of the former, adapted to apertures with an infinite planar geometry \cite{Devaney}.

Within the Kirchhoff–Helmholtz integral representation, the acoustic pressure wavefield in a solution domain bounded by a surface is expressed as an integral involving a linear combination of the pressure field and its outward normal derivative evaluated on that surface. However, because these two boundary quantities are linked by an additional homogeneous surface integral relation in a complementary domain, the radiated wavefield can equivalently be represented using only a single boundary quantity—namely, either the pressure or its normal derivative. Consequently, using these surface integral formulas, boundary conditions are imposed on the chosen solution space, in addition to enforcing Cauchy conditions that satisfy causality \cite{Devaney}. The choice of solution space and boundary conditions depends on the physics of the problem.

From a physical point of view, we consider two cases: 1) A vibrating piston is mounted in a \textit{rigid baffle}, where the normal derivative \( \frac{\partial p}{\partial \bn} \) of the pressure field vanishes everywhere on the baffle except on the acoustic aperture, which corresponds to the front face of the piston. In this configuration, a \textit{monopole} (or \textit{singlet}) source is defined in terms of the negative normal derivative of the pressure, \( -\frac{\partial p}{\partial \bn} \), or equivalently, \( \rho_0 \frac{\partial \bu}{\partial t} \cdot \bn \)~\cite{Li2}, where \( p \) denotes the scalar pressure, \( \rho_0 \) is the ambient density, \( \bu \) is the particle velocity vector, and \( \bn \) is the outward-pointing unit normal vector to the surface. This gives rise to the monopole integral representation, in which the Green’s function acts on a monopole source term restricted to the surface. 2) Under a \textit{soft-baffle} assumption, the pressure \( p \) is assumed to vanish on the baffle except on the acoustic aperture. In this case, the wavefield is represented using a surface-restricted \textit{dipole}\footnote{In the mathematical literature, the corresponding ``monopole'' and ``dipole'' integral formulas are typically expressed as ``single-layer'' and ``double-layer'' potential integral formulas, respectively.}(or \textit{doublet}) source, given by \( p \bn \). The soft-baffle representation involves the action of the normal derivative of the Green's function on \( p \), giving rise to the \textit{dipole} integral formula.

Interesting works have been conducted on the analytical calculation and numerical approximation of surface emitters and receivers, as well as their application to inverse problems \cite{Kostli,Finch,Xu,Burgholzer1,Burgholzer2,Cakoni,Stefanov,
Hristova,Ammari,Qian,Bao,Kunyanski,Dean-Ben,Huang,Haltmeier1,Belhachmi,Nguyen,Datchev,Ammari,Haltmeier2,Do,Uhlmann,Kuchment}, with some studies assuming smoothness of the corresponding surfaces. Some of these works, particularly in the context of photoacoustic tomography, have been discussed in the Introduction of~\cite{Javaherian10}.

For numerical approximations of the surface integral formulas, we employ a system of three coupled first-order wave equations that describe the propagation of acoustic waves in free space-time~\cite{Mast,Tabei}. Correspondingly, we define an \emph{extension operator} that maps a monopole source supported on an infinite plane to a scalar-valued mass source term, whose effective support lies within a narrow region surrounding the finite-sized support of the source on the plane. This mass source term is incorporated into the equation of continuity in free space-time domain.

In the case of the dipole integral formulation, the integrand involves the normal derivative of the Green's function evaluated on the aperture surface. Analytically, this normal derivative depends on the \emph{obliquity factor}, or equivalently, the \emph{solid angle}, which quantifies the angle subtended by an elemental surface area as observed from a given point in the domain. Accordingly, in the dipole case, the extension operator maps a vector-valued dipole source supported on a surface to a near-surface force source term, which is added to the equation of motion in the wave system.

Most importantly, building upon this extension operator, this work introduces a \emph{restriction operator} that maps the wavefield—solutions of the wave equation—to boundary data corresponding to the pressure (or its normal derivative) measured at the receiver surface. We then show for the practical cases of Dirichlet-type
boundary data (pressure) that the adjoint of a forward operator incorporating this trace mapping coincides---up to a constant factor---with a time-reversed formulation of the inwardly-directed dipole surface integral formula, evaluated on the receiver surface. The corresponding time-reversed system acts as a back-projection operator that inherently accounts for the frequency-dependent and analytically characterized angular sensitivity of finite-size receivers. Given that modeling acoustic apertures as omnidirectional point sources may lack sufficient accuracy—particularly in high-frequency regimes relevant to biomedical applications—the proposed adjoint operator offers a compelling alternative for incorporation into iterative frameworks, such as error minimization algorithms~\cite{Qian,Dean-Ben,Haltmeier1,Datchev,Arridge,Javaherian1}, to solve inverse problems involving receivers of finite size~\cite{Xu,Burgholzer1}.

 \textbf{Outline.} Section~2 introduces the wave equation in the time domain and explains how a unique solution can be obtained using homogeneous Cauchy conditions, which establish a causal relationship between the solution wavefield and a radiation source. A \textit{primary} solution, describing the propagated wavefield in terms of a radiation source with $d$-dimensional support (a volumetric region for $d=3$) in full space-time domain, is presented. Integral formulas expressing the wavefield in terms of a surface source are then derived. It will be demonstrated how imposing either \textit{Dirichlet} or \textit{Neumann} boundary conditions yields corresponding solutions to the overdetermined system arising from these integral representations.

Section~3 reformulates these surface integral formulas, upon which the 
boundary conditions are imposed, as the action of a causal Green's function 
on a source supported on an infinite plane.

Section~4 outlines the wave equation system underlying the derived time-domain analytical formulas. This section provides a detailed description of integral expressions that map the field or its normal derivative, supported on a $(d\!-\!1)$-dimensional infinite plane, to mass and force source terms supported in the free space-time domain, serving as inputs to the wave equation system.

In Section~5, building upon these integral formulations, we introduce the forward operator in the context of ultrasound tomography. In this section, we derive a \emph{restriction operator} that maps the wavefield—approximated in free space—onto the Neumann or Dirichlet boundary data supported on subsurfaces corresponding to the finite-sized receivers. Each of these subsurfaces is assumed to lie on a corresponding infinite plane. 

In Section~6, the adjoint of the resulting forward operator, incorporating the trace mappings, is derived using the integral representations developed in the preceding sections.

Section~7 explains the full-waveform approximation of the derived analytic integral formulas, discretized in time and on a regular grid.

Finally, Section~8 presents the numerical results obtained using this discretization process, while Section~9 discusses the broader significance of the study, particularly in the context of inverse problems.

\section{Wave Equation in the Time Domain }
\label{sec:wave-time}

\noindent
This section considers the propagation of acoustic waves from a real-valued space-and-time-varying source in an infinite, isotropic, and homogeneous medium in free space. Let $\bx = \left[x^1, \dots, x^d\right]^T$ denote a spatial position in $\mathbb{R}^d$ with $d \in \{2,3\}$ as the number of dimensions. The analysis presented here is performed for $d=3$ but holds for $d=2$ by replacing volumes with surfaces and surfaces with lines, and using the 2D Green's function. (For instance, a line source in a 3D medium produces acoustic waves that propagate as cylindrical waves, equivalent to an omnidirectional point source in a 2D medium.) 

\noindent
The real-valued wavefield satisfies the inhomogeneous wave equation, expressed as
\begin{align}  \label{eq:wave1}
   \bigg[\frac{1}{c^2}\frac{\partial^2}{\partial t^2} - \rho_0 \nabla \cdot \left( \frac{1}{\rho_0} \nabla \right)  \bigg] p(\bx, t)   = s(\bx,t).
\end{align}

\noindent
Here, the term on the right-hand side, $s$, is the forcing term, referred to here as the \textit{radiation source} (units: $\mathrm{kg} \, \mathrm{m}^{-\mathrm{d}} \, \mathrm{s}^{-2}$), and
is compactly supported in the spatio-temporal region 
\begin{align}  \label{eq:source-support}
    \Lambda_s = \{ \bx \in \nu_s \subset \mathbb{R}^d, \ t \in 
[ 0,T_s] \}.
\end{align}
Here, $\nu_s$ is a $d$-dimensional space (a volumetric region for $d=3$), and $[ 0, T_s]$ denotes the radiation time of the source. The radiation source $s$ is assumed to be square-integrable over $\Lambda_s$. 

\noindent
Additionally, $c$ represents the velocity of wave propagation in the medium (units: $\mathrm{ms}^{-1}$), and $\rho_0$ denotes the ambient density of the medium (units: $\mathrm{kg} \, \mathrm{m}^{-\mathrm{d}}$). The pressure wavefield $p$, the unknown parameter of the wave equation, has units of $\mathrm{kg} \, \mathrm{m}^{2-\mathrm{d}} \, \mathrm{s}^{-2}$ (or Pascal). 

\noindent
Assuming that \( \rho_0 \) is weakly heterogeneous, the wave equation~\eqref{eq:wave1} simplifies to its canonical form:

\begin{align}  \label{eq:wave2}
    \bigg[  \frac{1}{c^2} \frac{\partial^2 }{\partial t^2} - \nabla^2 \bigg] p(\bx, t) = s(\bx, t).
\end{align}

\subsection{Cauchy Conditions for Unique Solution}
Typically, the solution to the wave equation is nonunique. A unique solution is obtained by confining the wavefield to a particular solution that is causally related to the source, i.e., a wavefield $p$ that vanishes prior to the initial time $t=0$ of the source radiation. By imposing the Cauchy conditions

\begin{align}  \label{eq:cauchy1}
    p_+(\bx, t) \big|_{t=0} = 0, \quad 
    \frac{\partial p_+}{\partial t} (\bx, t) \big|_{t=0} = 0,
\end{align}
\noindent
a causal solution to the wave equation \eqref{eq:wave2} is obtained, where the subscript $+$ denotes causality. Solving the inhomogeneous wave equation \eqref{eq:wave2} with the causality conditions (Cauchy conditions \eqref{eq:cauchy1}) is referred to as the \textit{radiation problem}. 

By contrast, if the forcing term $s$ is set to zero, and the Cauchy conditions \eqref{eq:cauchy1} are replaced with arbitrary and inhomogeneous fields at $t=0$, the wave equation \eqref{eq:wave2} becomes an \textit{initial-value problem}.

\subsection{Green's Function Solution to the Wave Equation} \label{sec:gr}
Consider the wave equation \eqref{eq:wave2} as a radiation problem for a particular choice of source, $s(\bx, t) = \delta(\bx - \bx') \delta(t - t')$, where $\delta(\cdot)$ is the Dirac delta distribution, where $\bx'$ and $t'$ are free parameters in the space and time domains, respectively. Here, $\bx$ and $t$ are fixed, but arbitrary parameters in a chosen spatio-temporal region $\Lambda = \{\bx \in \nu, \ t \in (t_0, t_1)\}$. In an infinite free space-time, the Green's function solution to the wave equation satisfies

\begin{align}  \label{eq:greens}
    \bigg[ \frac{1}{c^2} \frac{\partial^2}{\partial t^2} - \nabla^2 \bigg] g(\bx - \bx', t - t') = \delta(\bx - \bx') \delta(t - t').
\end{align}

\noindent
For brevity, the Green's function will henceforth be denoted as
\begin{align} \label{eq:recipro}
    g(\bx_{\fd}, t_{\fd}) := g(\bx - \bx',\, t - t'),
\end{align}
where $\bx_{\fd} = \bx - \bx'$ and $t_{\fd} = t - t'$.
 Similar to the wave equation \eqref{eq:wave2}, a unique solution to Eq.~\eqref{eq:greens} is obtained by assuming a causality condition for the Green's function, i.e., $g_+(\bx_{\fd}, t_{\fd}) = 0$ for $t_{\fd} < 0$. For $d=3$, the causal Green's function satisfies in the free space-time

\begin{align} \label{eq:greens-3d-time}
    g_+(\bx_{\fd}, t_{\fd}) = \frac{1}{4 \pi} \frac{\delta(t_{\fd} - \frac{x_{\fd}}{c})}{x_{\fd}}.
\end{align}
\noindent
Rewriting the wave equations \eqref{eq:wave2} and \eqref{eq:greens} in the forms

\begin{align}  \label{eq:wave4}
    \bigg[ \frac{1}{c^2} \frac{\partial^2}{\partial {t'}^2} - \nabla_{\bx'}^2 \bigg] p_{+}(\bx', t') = s(\bx', t'),
\end{align}

\noindent
and 

\begin{align}  \label{eq:greens3}
    \bigg[ \frac{1}{c^2} \frac{\partial^2}{\partial {t'}^2} - \nabla_{\bx'}^2 \bigg] g(\bx - \bx', t - t') = \delta(\bx - \bx') \delta(t - t').
\end{align}
\noindent

Now, multiplying Eq.~\eqref{eq:wave4} by $g(\bx_{\fd}, t_{\fd})$ and Eq.~\eqref{eq:greens3} by $p_+(\bx', t')$, then subtracting the modified Eq.~\eqref{eq:greens3} from the modified Eq.~\eqref{eq:wave4}, yields \cite{Devaney}:

\begin{align} \label{eq:gr-th1}
    \frac{1}{c^2} \bigg[g  \big[ \frac{\partial^2 p_+}{\partial {t'}^2} \big]
    - \big[ \frac{\partial^2 g}{\partial {t'}^2} \big] p_+ \bigg] 
    - \bigg[g \big[ \nabla^2 p_+\big] - \big[ \nabla^2 g \big] \, p_+ \bigg] = g s - \delta(\bx_{\fd}) \delta(t_{\fd}) \, p_+.
\end{align}

Integrating Eq.~\eqref{eq:gr-th1} on both sides over the spatio-temporal solution region \( \Lambda \), taking the temporal integral of the first term in \eqref{eq:gr-th1}, and applying the divergence theorem to the second term yield \cite{Devaney}:
\begin{align} \label{eq:gr-th2}
\begin{split}
    & \frac{1}{c^2} \int_\nu d\bx' \ 
    \bigg[g \big[ \frac{\partial p_+}{\partial t'} \big] - \big[ \frac{\partial g}{\partial t'} \big] p_+ \bigg] 
    \bigg|_{t' = t_0}^{t_1} 
    - \int_{t_0}^{t_1} \! dt' \int_{\partial \nu} \! dS' \ 
    \bigg[  g \big[ \frac{\partial p_+}{\partial \bn} \big]  - \big[ \frac{\partial g}{\partial \bn} \big] p_+  \bigg] = \\
    & \int_{t_0}^{t_1} \! dt' \int_\nu \! d\bx' \ g s - 
    \begin{cases}
        p_+, & \text{if } \bx \in \nu, t \in (t_0, t_1), \\
        0, & \text{otherwise},
    \end{cases}
\end{split}
\end{align}
where $\bn$ denotes the outward unit normal to the surface 
$\partial \nu$, pointing from the region $\nu$ into its complement.

\subsubsection{Primary Solution} \label{sec:gr-pri}
A primary solution for all space and all time can be obtained by setting \( g \) the causal Green's function, \( g_+ \), and extending the spatio-temporal set \( \Lambda \) to infinity. Accordingly, the limits \( t_0 \to -\infty \), and \( t_1 \to \infty \) are taken, and the region \( \nu \) extends to infinity in radius, \( R_{\infty} \).

For the first term on the left-hand side of Eq.~\eqref{eq:gr-th2}, the assumption of causality ensures that $p_+$ vanishes at $t'= -\infty$, and $g$ vanishes at $t' = +\infty$. As a result, this term drops out. For the latter, we have used the fact that a causal Green's function satisfies $g(\bx_{\fd}, t_{\fd}) = 0$ unless $x_{\fd} = c t_{\fd}$. For the second term, due to the same reasoning, the contribution from the surface at infinite radius (\( R_\infty \)) vanishes~\cite{Devaney}. Consequently, the second term on the left-hand side of Eq.~\eqref{eq:gr-th2} also vanishes.
\noindent
Therefore, the primary solution for $p_+$ at any pair of $\bx$ and $t$ lying in the domain $\Lambda$ satisfies

\begin{align} \label{eq:gr-prime}
    p_+(\bx, t) = \int_{\mathbb{R}} dt' \int_{\mathbb{R}^d} d\bx' \, g_+(\bx_{\fd}, t_{\fd}) s(\bx', t'),
\end{align}
where $s \in C_0^{\infty}(\Lambda_s ) $, and $\Lambda$ is, in this case, set to be the full space--time domain.

\subsubsection{Kirchhoff-Helmholtz Solution} \label{sec:gr-kh}
This section describes the \textit{Kirchhoff--Helmholtz solution} to the wave equation. While the primary solution in Eq.~\eqref{eq:gr-prime} directly maps the forcing term $s$ to the wavefield solution $p_+$ via an integral over the volumetric support region $\nu_s$, many practical problems instead involve sources distributed over a bounding surface.

Consider a region \( \nu \subset \mathbb{R}^d \) that contains the source support \( \nu_s \) and is bounded by a closed surface \( \partial \nu \). We now confine the solution space to the region \( \nu^+ \), which lies outside \( \partial \nu \), i.e., 
\[
\nu^+ := \mathbb{R}^d \setminus \left( \nu \cup \partial \nu \right).
\]
By assuming causality for the solution \( p_+ \) in Eq.~\eqref{eq:wave4} over the time interval \( (t_0, t_1) \), where \( t_0 \to -\infty \) and \( t_1 \to +\infty \), and by applying the causality condition for the Green's function \( g_+ \) in Eq.~\eqref{eq:greens3}, the solution procedure remains identical to that described in Section~\ref{sec:gr-pri}. However, the integral in Eq.~\eqref{eq:gr-th2} is now evaluated over the external region \( \nu^+ \) rather than \( \nu \).

As before, the first term on the left-hand side of Eq.~\eqref{eq:gr-th2} vanishes due to the causality of $p_+$ and $g_+$, and the contribution from the surface at infinite radius also vanishes, as discussed in Section~\ref{sec:gr-pri}. Consequently, the field $p_+$ in the solution space $\nu^+$ satisfies the following integral equation \cite{Devaney}:
\begin{align}  \label{eq:kh}
    \int_{-\infty}^{\infty} dt' \int_{\partial \nu} dS' \, 
    \bigg[ 
       \big[   \frac{\partial g_+}{\partial \bn} \big] p_+ 
        - g_+  \big[ \frac{\partial p_+}{\partial \bn} \big]
    \bigg] =
    \begin{cases}
        p_+(\bx, t) & \bx \in \nu^+\\
        0 & \bx \in \nu,
    \end{cases}
\end{align}
where $dS' := dS(\bx')$ denotes the elemental surface area at position $\bx'$, and $\bn$ is the unit normal vector to the surface $\partial \nu$, oriented outward from the interior volume $\nu$ toward the exterior region $\nu^+$.

In Eq.~\eqref{eq:kh}, the formula for the exterior volume is known as the \textit{first Helmholtz identity}. It expresses the wavefield \( p_+ \) outside the surface \( \partial \nu \) in terms of the field and its normal derivative on \( \partial \nu \). In the corresponding formula for the interior volume, referred to as the \textit{second Helmholtz identity}, the right-hand side vanishes, establishing a dependence between the field and its normal derivative on \( \partial \nu \). This dependence allows the integrand in Eq.~\eqref{eq:kh} to be reformulated solely in terms of either the wavefield or its normal derivative. These conditions are enforced by imposing \textit{Dirichlet} or \textit{Neumann} boundary conditions, leading to the \textit{dipole} and \textit{monopole} integral formulas, respectively.

\subsubsection{Rayleigh-Sommerfeld solution} \label{sec:ray-som}

The \textit{Rayleigh--Sommerfeld} integral formula arises from solving the boundary-value problem for the Helmholtz wave equation under the assumption of a source mounted on the bounding plane of a half-space~\cite{Devaney}. To derive this formula, the solution space is defined as a half-space bounded by an infinite plane, denoted by \( \partial \nu_{\text{plane}} \), and a hemisphere of infinite radius, denoted by \( R_{\infty} \), which bounds the half-space. For simplicity, the bounding plane is taken as \( x^d = 0 \), and the solution is sought in the half-space \( x^d > 0 \).

As discussed previously, the associated \textit{Kirchhoff-Helmholtz} integral formula is overdetermined. Therefore, boundary conditions of the \textit{Dirichlet} or \textit{Neumann} types are imposed on \( \partial \nu \). A common method for solving this boundary-value problem is the \textit{method of images} \cite{Devaney}. Here, \( \bx \) and \( \bx' \) represent the positions of the general field point and the source point in $\mathbb{R}^d$, respectively, with both assumed to lie in the half-space \( \nu^+ = \{\bx : x^d > 0\} \). 

For the source point \( \bx' = \big[ x'^1, \ldots, x'^d \big]^T \), a mirror-image source point \( \tilde{\bx}' = \big[ x'^1, \ldots, -x'^d \big]^T \) is introduced to ensure that \( \bx \) and \( \tilde{\bx}' \) lie on opposite sides of the plane, thereby satisfying \( \delta(\bx - \tilde{\bx}') = 0 \). Correspondingly, the Helmholtz equation is defined with an augmented forcing term:
\begin{align} \label{eq:holm-med}
   \left[ \frac{1}{c^2} \frac{\partial^2}{\partial t'^2} -  \nabla_{\bx}^2 \right] g_+^N(\bx - \bx', t_{\fd}) = \left[ \delta(\bx - \bx') + \delta(\bx - \tilde{\bx}') \right] \delta(t_{\fd}),
\end{align}
where the augmented Green's function \( g_+^N \) is given by:
\begin{align} \label{eq:neu}
    g_+^N(\bx - \bx', t_{\fd}) = g_+(\bx - \bx', t_{\fd}) + g_+(\bx - \tilde{\bx}', t_{\fd} ).
\end{align}

When the source point \( \bx' \) approaches the plane \( x^d = 0 \), or when the field point \( \bx \) itself lies on or near this plane, we have \( |\bx - \tilde{\bx}'| = |\bx - \bx'| \). Consequently, the Green's function satisfies:
\begin{align} \label{eq:neu1}
    g_+^N(\bx_{\fd}, t_{\fd}) = 2 g_+(\bx_{\fd}, t_{\fd}),
\end{align}
along with the Neumann boundary condition:
\begin{align}\label{eq:neu-b}
    \frac{\partial g_+^N}{\partial \bn}(\bx_{\fd}, t_{\fd}) = 0,
\end{align}
on \( \partial \nu = \{ \bx : x^d = 0 \} \).

\noindent
Since the Green's function satisfies the Neumann boundary condition on the infinite plane, this formulation yields the \textit{Rayleigh–Sommerfeld} monopole integral representation:
\begin{align} \label{eq:gr-rsn}
    p_+^N(\bx, t) = 
    -2 \int_{-\infty}^{\infty} dt' \int_{\partial \nu} dS' \ g_+(\bx_{\fd}, t_{\fd}) \left[ \frac{\partial p_+}{\partial \bn} (\bx', t')\right],
\end{align}
where \( \bx \in \nu^+ \), and \(\bn\) is the unit normal vector to the plane \(x^d = 0\), directed into the half-space \(\nu^+\).

\noindent
Alternatively, changing the sign of \( \delta(\bx - \tilde{\bx}') \) in the right-hand side of Eq.~\eqref{eq:holm-med} yields:
\begin{align} \label{eq:dir1}
    \frac{\partial g_+^D}{\partial \bn}(\bx_{\fd}, t_{\fd}) = 2 \frac{\partial g_+}{\partial \bn}(\bx_{\fd}, t_{\fd}),
\end{align}
with the Dirichlet boundary condition:
\begin{align}\label{eq:dir-b}
    g_+^D(\bx_{\fd}, t_{\fd}) = 0.
\end{align}

\noindent
This leads to the \textit{Rayleigh–Sommerfeld} dipole integral representation:
\begin{align} \label{eq:gr-rsd1}
    p_+^D(\bx, t) = 2 \int_{-\infty}^{\infty} dt' \int_{\partial \nu} dS' \ 
    \left[ \frac{\partial g_+}{\partial \bn} (\bx_{\fd}, t_{\fd})\right] p_+(\bx', t').
\end{align}


Equation~\eqref{eq:gr-rsn} corresponds to the \textit{rigid-baffle} condition, in which the normal derivative of the wavefield vanishes on the baffle surface, except on the acoustic aperture. Similarly, Eq.~\eqref{eq:gr-rsd1} corresponds to the \textit{soft-baffle} condition, where the wavefield itself vanishes everywhere on the baffle except on the aperture.

\section{Monopole and Dipole Formulas in Terms of Actions of the Causal Green's Function on the Surface Source} \label{sec:sum}

\noindent
This section derives the wavefield explicitly in terms of the action of the causal Green's function on sources supported on the infinite hyperplane bounding the half-space.

\subsection{Monopole Formula}

As previously discussed, the monopole integral formula is derived under the assumption of a vibrating piston mounted on a rigid baffle. This configuration enforces a vanishing normal derivative of the wavefield on the baffle, except on the front face of the vibrating piston. Consequently, a \textit{Neumann} Green's function, denoted by \( g_+^N \), is used. This Green's function satisfies a homogeneous Neumann boundary condition on the surface \( \partial \nu \). 

\noindent
By assuming that the surface is an infinite plane bounding a half-space, the surface integral formula reduces to a dependence solely on the normal derivative of the pressure field, namely
\(
-\,\frac{\partial p_+}{\partial \bn},
\)
which represents a monopole source. Accordingly, using the Neumann Green's function in Eq.~\eqref{eq:neu1}, the time-domain Rayleigh–Sommerfeld integral equation can be written in terms of the action of the causal Green’s function on the monopole source as

\begin{align} \label{eq:gr-monopole}
    p_+^N(\bx, t) = 2 \int_{-\infty}^{\infty} dt' \int_{\partial \nu} dS' \ g_+(\bx_{\fd}, t_{\fd}) \left[ \rho_0 \frac{\partial u^{\bn}}{\partial t} (\bx', t') \right],
\end{align}
where \( \bx \in \nu^+ \), and the identity
\begin{align}
    \frac{\partial p}{\partial \bn} = - \rho_0 \frac{\partial u^{\bn}}{\partial t},
\end{align}
with \( u^{\bn} = \bu \cdot \bn \), has been used. Here, \( \rho_0 \) is the ambient mass density, and \( u^{\bn} \) is the normal component of the velocity vector \( \bu \) on the surface.

\subsection{Dipole formula} \label{sec:dipole} 
The dipole formula is derived under the assumption of a soft baffle, on which the pressure vanishes everywhere except the front face of the vibrating piston. Correspondingly, using a causal Green's function $g_+^D$, which satisfies a homogeneous \textit{Dirichlet} boundary condition over $\partial \nu$—an infinite plane bounding the half-space—the integral formula depends solely on $p_+ \bn$, representing a dipole source. For analytic (or ray-based) methods used to approximate the dipole integral formula, it is convenient to reformulate the formula in terms of weighted actions of the causal Green's function. Accordingly, using the Dirichlet Green's function in Eq.~\eqref{eq:dir1}, the dipole variant of the \textit{Rayleigh-Sommerfeld} formula can be expressed in terms of action of the spatial derivative of the causal Green's function acting on the dipole source as 
\begin{align} \label{eq:gr-dipole1}
\begin{split}
    p_+^D(\bx, t) &= 2  \int_{-\infty}^{\infty} dt' \ \int_{\partial \nu} dS' \   \nabla_{\bx'}  g(\bx_{\fd}, t_{\fd}) \cdot \Big[ p(\bx', t') \bn \Big] \\
    & = 2 \int_{-\infty}^{\infty} dt' \ \int_{\partial \nu} dS' \  \Big[ \bn \cdot \frac{\bx_{\fd}}{x_{\fd}}\Big]  \ g_+(\bx_{\fd}, t_{\fd}) \frac{1}{c}\Big[ \frac{\partial}{\partial t'} + \frac{1}{t_{\fd}} \Big] p(\bx', t'), 
\end{split}
\end{align}
where $\bx \in \nu^+$ and $\bx_{\fd} = \bx - \bx'$ is the distance vector. In the second line of this formula, $\bn \cdot \bx_{\fd}/x_{\fd}$ is the \textit{obliquity factor}, which weights the actions of the causal Green's function on a monopole-like source, decomposed into \textit{far-field} and \textit{near-field} terms. Applying the \textit{far-field approximation}, where $ \frac{\partial p}{\partial t'} \gg \frac{p}{t_{\fd}} $, which is equivalent to $k x_{\fd} \gg 1$ in the frequency domain, and valid in regions sufficiently far from the aperture, allows neglecting the $p/t_{\fd}$ term. Thus, the \textit{far-field approximation} of Eq.~\eqref{eq:gr-dipole1} becomes:
\begin{align} \label{eq:gr-rsd}
    p_+^D(\bx, t) \approx   2  \int_{-\infty}^{\infty} dt' \ \int_{\partial \nu} dS' \     \Big[ \bn \cdot \frac{\bx_{\fd}}{x_{\fd}} \Big] \ g_+(\bx_{\fd}, t_{\fd}) \ \bigg[  \frac{1}{c}  \frac{\partial p}{\partial t'}  (\bx', t')  \bigg].
\end{align}

\begin{remark}
Using ray-based methods, it is computationally more efficient to express the integral formula \eqref{eq:gr-dipole1} in terms of a solid angle element $d \Omega_{\bx}\big(S(\bx')\big)$, defined as the angle subtended by an infinitesimal area $dS$ corresponding to the point $\bx'$ on the surface of the aperture, as seen from any arbitrary field point $\bx \in \nu^+$. The solid angle relates to the obliquity factor through the formula:
\begin{align} \label{eq:solid}
    d  \Omega_{\bx}\big(S(\bx')\big) = \frac{dS'}{x_{\fd}^2} \ \Big[ \bn \cdot \frac{\bx_{\fd}}{x_{\fd}} \Big].
\end{align}
\end{remark}

In the next section, we demonstrate how a scalar-valued mass source or a vector-valued force source can be defined and incorporated into the equations of continuity and motion, respectively. This formulation ensures that the resulting system of wave equations approximates the monopole formula \eqref{eq:gr-monopole} and the dipole integral formula \eqref{eq:gr-dipole1}, respectively.

\section{Full-Waveform Approximation of the Wave Equation in the Time Domain}
\label{sec:approx}

\noindent
This section presents a full-waveform approximation of the acoustic wave equation~\eqref{eq:wave1}. In particular, it formulates a wave-equation system for approximating the time-domain wavefield induced by a source distribution supported on an infinite plane bounding a half-space. The wavefield is represented based on the monopole and dipole integral representations given in Eqs.~\eqref{eq:gr-monopole} and~\eqref{eq:gr-dipole1}.

\subsection{System of Wave Equations in Terms of Source \( s \)}
\label{sec:semi-num-s}

\noindent
To better understand these integral representations, we begin by replacing the wave equation~\eqref{eq:wave1} with a coupled first-order linear system:
\begin{align}  
\label{eq:stepping1}
\begin{split}
    & \frac{\partial}{\partial t} \bu (\bx, t) = - \frac{1}{\rho_0} \nabla p(\bx, t),  \\
    & \frac{\partial}{\partial t} \rho(\bx,t) = - \rho_0 \nabla \cdot \bu (\bx, t) + s_m(\bx, t), \\
    & p(\bx, t) = c^2 \rho(\bx,t),
\end{split}
\end{align}
where \( s_m \in C_0^{\infty}(\Lambda_s) \) denotes the mass source, assumed to satisfy
\begin{align} \label{eq:mass-source}
s(\bx,t) = \frac{\partial }{\partial t} s_m(\bx,t),
\end{align}
for some source \( s \in C_0^{\infty}(\Lambda_s)  \). Recall that \( \Lambda_s \) denotes the support of the source distribution, as defined in Eq.~\eqref{eq:source-support}. In the wave propagation model, this distribution is extended by zero to the free space.

\subsection{System of Wave Equations in Terms of Regularized Source \( \mathcal{S} \)}
\label{sec:semi-num}

The wave system~\eqref{eq:stepping1} can be used directly to compute the time-domain \textit{primary solution} given in Eq.~\eqref{eq:gr-prime}, which is valid over the full space-time domain. However, the source \( s \in C_0^{\infty}(\Lambda_s) \) is not physically accessible. Instead, the wavefield should be expressed using an equivalent source supported on an infinite plane that bounds the half-space, based on the monopole and dipole integral representations in Eqs.~\eqref{eq:gr-monopole} and~\eqref{eq:gr-dipole1}.

\noindent
Accordingly, we define $\mathcal{S}$ as a regularized source, whose effective support lies within a narrow region surrounding the surface and decays smoothly away from it. To this end, we rewrite the wave equation~\eqref{eq:wave1} by moving the second term on the left-hand side to the right-hand side, yielding
\begin{align}  
\label{eq:wave3}
  \frac{1}{c^2} \frac{\partial^2 p (\bx, t)}{\partial t^2}
  = \rho_0 \nabla \cdot \left( \frac{1}{\rho_0} \nabla p(\bx, t) \right) 
  + \mathcal{S}(\bx, t),
\end{align}
where \( \mathcal{S} \in C_0^{\infty} \big( \mathbb{R}^d \times \mathbb{R}^+ \big) \) denotes a source term defined over the full space-time domain that remains square-integrable in both space and time.

\noindent
A coupled first-order linear reformulation of the wave equation~\eqref{eq:wave3} leads to the following system \cite{Tabei}:
\begin{align}  
\label{eq:stepping3}
\begin{split}
    & \frac{\partial}{\partial t} \bu (\bx, t) = - \frac{1}{\rho_0} \nabla p(\bx, t) + \boldsymbol{\mathcal{S}}_f(\bx, t), \\
    & \frac{\partial}{\partial t} \rho(\bx,t) = - \rho_0 \nabla \cdot \bu (\bx, t) + \mathcal{S}_m(\bx, t), \\
    & p(\bx, t) = c^2 \rho(\bx,t).
\end{split}
\end{align}

This system solves the regularized wave equation \eqref{eq:wave3} in free space-time domain. Accordingly, the first line of the system represents the equation of \textit{motion}, where $\boldsymbol{\mathcal{S}}_f \in C_0^{\infty}( \mathbb{R}^d \times \mathbb{R}^+ )$ denotes a vector-valued force source. The second line corresponds to the equation of \textit{continuity}, with $\mathcal{S}_m \in C_0^{\infty}( \mathbb{R}^d \times \mathbb{R}^+ )$ representing a scalar-valued mass source, which satisfies
\begin{align} \label{eq:mass-source-regularized}
    \mathcal{S}(\bx, t) = \frac{\partial}{\partial t} \mathcal{S}_m(\bx, t).
\end{align}
\noindent
All source terms will later be expressed as integral formulas involving source distributions compactly supported on a surface. To achieve this, a regularized Dirac delta distribution is introduced.

\begin{definition} \label{def:1}
 Let $x^{\zeta}$ denote the coordinate of $\bx \in \mathbb{R}^d$ in the Cartesian direction $\zeta \in \{1,\dots,d\}$. The regularized Dirac delta distribution is defined as
\begin{align} \label{eq:delta-regularized}
\delta_b(\bx - \bx') 
= \prod_{\zeta=1}^d \frac{1}{b} \;
\mathrm{sinc}\!\left( \frac{\pi (x^\zeta - x'^\zeta)}{b} \right),
\end{align}
where $b > 0$ is a bandwidth parameter that controls the spread of the regularization (units: $\mathrm{m}^{-\mathrm{d}}$). Here, the sinc function is defined by
\[
\sinc(f) = \frac{\sin(f)}{f}, \quad \text{with} \quad \sinc(0) = 1.
\]
The function \( \delta_b \) defined in~\eqref{eq:delta-regularized} converges to the Dirac delta distribution in the limit \( b \to 0^+ \), satisfying the sifting property
\begin{align} \label{eq:sifting}
\lim_{b \to 0^+} \int_{\mathbb{R}^d} d\bx' \, \delta_b(\bx - \bx') \, f(\bx') = f(\bx).
\end{align}
\end{definition}

\begin{remark}
Although the forward problem is posed in the continuous domain, our ultimate goal is to discretize the formulas on a regular grid. Since such discretization inherently suppresses high-frequency components, we work with infinitely differentiable fields and regularized Dirac delta distributions (see Eq.~\eqref{eq:delta-regularized}), with the bandwidth parameter \( b \to 0^+ \). In this way, the fields remain smooth but converge, in the distributional sense, to their counterparts that are discontinuous on the infinite plane. This allows us to retain essential analytical properties, such as the sifting property of the Dirac delta distribution, while correctly handling discontinuities and singularities.
\end{remark}

\subsection{Definition of the Regularized Source \( \mathcal{S} \) in Terms of the Surface Integral Formulas}

This section employs the system of coupled first-order wave equations \eqref{eq:stepping3} to describe the wavefield in terms of surface sources via the monopole and dipole integral formulas. In particular, it shows how the inclusion of a mass source $\mathcal{S}_m$ in the continuity equation within the wave system \eqref{eq:stepping3} enables approximation of the monopole integral formula \eqref{eq:gr-monopole}. Additionally, it demonstrates how incorporating a vector-valued force source $\boldsymbol{\mathcal{S}}_f$ in the equation of motion allows the system to represent the dipole integral formula \eqref{eq:gr-dipole1}.

\subsubsection{Wavefield Representation in Terms of \( \rho_0 \frac{\partial \bu}{\partial t} \cdot \bn \) Compactly Supported on an Infinite Plane}

This subsection outlines how the wave system~\eqref{eq:stepping3} represents the wavefield in terms of a monopole source \(  \rho_0 \frac{\partial \bu}{\partial t} \cdot \bn \), compactly supported on an infinite plane. To achieve this, a regularized mass source is constructed from the surface-supported source as follows:
\begin{align} \label{eq:mass-u}
    \mathcal{S}_{(m,u^{\bn})}(\bx, t) = 2
    \int_{\partial \nu} dS(\bx_s) \ \delta_b(\bx - \bx_s) \Big[ \rho_0(\bx_s)  \, \bu(\bx_s, t) \cdot \bn \Big],
\end{align}
where \( \bn \) denotes the unit normal to the infinite plane \( \partial \nu \), directed into the half-space \( \nu^+ \). Under this approximation, the monopole integral representation in Eq.~\eqref{eq:gr-monopole} is reformulated as Eq.~\eqref{eq:gr-prime}, where the original source term \( s \) is replaced by the regularized near-surface source \( \mathcal{S} \).

Accordingly, the system \eqref{eq:stepping3}, with a regularized mass source defined by Eq. \eqref{eq:mass-u}, solves the regularized wave equation:
\begin{align}  \label{eq:wave-monopole}
\begin{split}
   \bigg[\frac{1}{c^2}\frac{\partial^2}{\partial t^2} - \rho_0 \nabla \cdot \left(\frac{1}{\rho_0} \nabla \right) \bigg] p(\bx, t)  & =  2 \  \int_{\partial \nu} dS(\bx_s) \  \delta_b(\bx-\bx_s)   \Big[ \rho_0(\bx_s)   \frac{\partial \bu}{\partial t} (\bx_s, t) \cdot \bn \Big]\\
   & = 2  \  \int_{\partial \nu} dS(\bx_s) \  \delta_b(\bx-\bx_s)   \Big[- \frac{\partial p}{ \partial  \bn}(\bx_s, t) \Big].
   \end{split}
\end{align}

\subsubsection{Wavefield Representation in Terms of a Monopole-like Source \( p \) Compactly Supported on an Infinite Plane} 
\label{sec:doublet}

In some studies, the source expressed in terms of \( p \), supported on a surface, is treated as a monopole source and approximated using a regularized scalar-valued mass source definition under the following assumptions:
\begin{enumerate}
    \item The far-field approximation—where \( \frac{\partial p}{\partial t'} \gg \frac{p}{t_{\fd}} \)—leading to the reduced dipole formula~\eqref{eq:gr-rsd}.
    \item An additional assumption \( \bn \cdot \bx_{\fd}/x_{\fd} \approx 1 \) in the reduced formula~\eqref{eq:gr-rsd}, which effectively treats the acoustic aperture as omnidirectional.
\end{enumerate}
Under these two conditions, the wavefield is approximated using the system~\eqref{eq:stepping3}, together with a regularized mass source:
\begin{align} \label{eq:mass-p2}
    \mathcal{S}_{(m,p)}(\bx, t)  
    = 2 \int_{\partial \nu} dS(\bx_s) \ \delta_b(\bx - \bx_s) \left[ \frac{1}{c(\bx_s)} \, p(\bx_s, t) \right].
\end{align}

\noindent
It is important to note that these assumptions may break down when \( \bx_s \) lies on a finite-sized surface, and for this reason, the approach is generally not recommended. Nevertheless, we employ it here as a benchmark to compare with the the dipole integral formula in Eq.~\eqref{eq:gr-dipole1}, demonstrating how the force-source definition introduced below effectively captures directionality and near-field effects in the wavefield approximation.

\subsubsection{Wavefield Representation in Terms of a Dipole Source \( p \bn \) Compactly Supported on an Infinite Plane}
\label{sec:force-source}

This subsection explains how the system of coupled first-order wave equations~\eqref{eq:stepping3} realizes the time-domain dipole integral formula~\eqref{eq:gr-dipole1}. To this end, a vector-valued regularized force source is constructed in terms of the dipole source \( p \bn \), supported on an infinite plane. This near-surface supported force source is defined as:
\begin{align} \label{eq:force}
\boldsymbol{\mathcal{S}}_f(\bx, t) =  \frac{2}{\rho_0(\bx)} \ 
\int_{\partial \nu} dS(\bx_s) \, \delta_b(\bx - \bx_s) \, \left[ \, p(\bx_s, t)    \, \bn \right],
\end{align}
where \( \bn \) is the outward-pointing unit normal vector on the surface \( \partial \nu \), directed into the region \( \nu^+ \).

With this definition of the force source, it can be shown via straightforward algebra that the system~\eqref{eq:stepping3} is equivalent to the regularized second-order wave equation:
\begin{align} \label{eq:wave-dipole}
\left[\frac{1}{c^2} \frac{\partial^2}{\partial t^2} - \rho_0 \nabla \cdot \left( \frac{1}{\rho_0} \nabla \right) \right] p(\bx, t) = -\rho_0 \nabla \cdot \boldsymbol{\mathcal{S}}_f(\bx,t) =
- 2 \ \int_{\partial \nu} dS(\bx_s) \, \nabla \delta_b(\bx - \bx_s) \cdot \Big[ p(\bx_s, t) \, \bn \Big] .
\end{align}

The wave system~\eqref{eq:stepping3}, with the force source \( \boldsymbol{\mathcal{S}}_f \) defined in Eq.~\eqref{eq:force}, plays a central role in the analysis of inverse problems. In the next section, we show that when boundary data are of Dirichlet type—i.e., given in the form of pressure measurements—the adjoint of the wave equation~\eqref{eq:wave1} corresponds to a time-reversed version of Eq.~\eqref{eq:wave-dipole}, with the source term appropriately time-reversed as well.

\section{Forward Problem Formulation in the Ultrasound Tomography Context} \label{sec:forward}

Building on the integral representations derived in the previous sections, we now formulate the forward problem. A general and representative example of this forward problem arises in ultrasound tomography, which can then be specialized to other applications such as therapeutic ultrasound optimization and photoacoustic tomography. In this study, we focus on the formulation of the forward and adjoint (or time-reversal\footnote{The notion of an \emph{adjoint} operator frequently appears in the context of iterative minimization of an objective function; by contrast, the terms \emph{time-reversal}, \emph{backprojection}, or \emph{inverse} operator all refer to the same concept and are more closely associated with single-step approaches, although embedding these latter operators within iterative minimization frameworks, such as \emph{Neumann-series} iterations, has been reported \cite{Qian}.
}) operators in this context, with the aim of embedding them into either an iterative inverse-problem framework or a single-step backprojection operator. Building on these derived operators and their equivalent full-waveform approximations, the associated minimization problems can be solved using first-order or second-order optimization methods \cite{Arridge,Haltmeier1,Javaherian1,Javaherian2}, or within Neumann-series frameworks \cite{Qian}. Exploration of the latter is left to future studies.

\subsection{Problem Setting in Ultrasound Tomography}  \label{sec:forward-def}

It was shown that the wavefield in a half space can be represented via an integral formulation involving the field or its normal derivative evaluated on a bounding surface \( \partial \nu \), where \( \partial \nu \) is taken to be an infinite plane \cite{Devaney}.

In this section, we assume that \( \Omega \subset \mathbb{R}^d \) is an open domain with boundary \( \partial \Omega \), on which both the emitted and measured data are prescribed. We decompose the full space into the interior domain \( \Omega^- := \Omega \) and the exterior domain \( \Omega^+ := \mathbb{R}^d \setminus \big( \Omega^- \cup \partial \Omega \big) \).

Accordingly, we index emitters and receivers by \( e \) and \( r \), respectively. Each emitter or receiver is treated as an infinite plane with support confined to a finite-sized subsurface lying on \( \partial \Omega \), using the infinite-plane assumption introduced in Section~\ref{sec:ray-som}. The infinite plane associated with each emitter \( e \) or receiver \( r \) is denoted by \( \partial \nu_e \) or \( \partial \nu_r \), respectively. Each infinite plane splits the domain into a half-space \( \nu_{\{e,r\}}^+ \), which contains \( \Omega^- \), and a complementary half-space \( \nu_{\{e,r\}}^- := \mathbb{R}^d \setminus \big( \nu_{\{e,r\}}^+ \cup \partial \nu_{\{e,r\}} \big) \). Consequently, the emitted and measured boundary data are compactly supported on these finite-sized subsurfaces \( \partial \nu_{\{e,r\}} \cap \partial \Omega \), which are often disk-shaped in practical implementations.

Each emitter--receiver pair produces a distinct time trace. Spherical or cylindrical acquisition geometries for \( \partial \Omega \) are often realized by assembling multiple such transducers. Accordingly, in this setting, the forward operator can be represented as one of the following boundary-to-boundary maps: Dirichlet-to-Neumann, Dirichlet-to-Dirichlet, Neumann-to-Neumann, or Neumann-to-Dirichlet. Among these, the \emph{to-Dirichlet} forward maps are most frequently encountered in practical applications and constitute the primary focus of analysis in this study. (We will also consider the \emph{to-Neumann} maps to preserve generality.)

\subsection{Forward Operator in Ultrasound Tomography}

Since the unknown parameters in inverse problems—such as those arising in ultrasound tomography—are typically defined over the full space-time domain, with constraints imposed in the exterior region, the forward operator is often reduced to a mapping
\begin{align}
    \begin{split}
        \mathcal{A}: C_0^{\infty} \left( \mathbb{R}^d \times \mathbb{R}^+ \right) &\rightarrow C_0^{\infty} \left( \partial \Omega \times [0, T] \right), \\
        \mathcal{A}_r \left[ \mathcal{S}_e \right] &= y_{r,e},
    \end{split}
\end{align}
where \( \mathcal{S}_e \in C_0^{\infty}(\mathbb{R}^d \times \mathbb{R}^+) \) denotes the regularized source term given by the right-hand side of the wave equations~\eqref{eq:wave-monopole} and~\eqref{eq:wave-dipole}, with integration taken over the emitting surface \( \partial \nu := \partial \nu_e \cap \partial \Omega \). Although the source \( \mathcal{S}_e \) is physically active only over a finite time interval, we formally extend its definition to all of space-time by assuming \( \mathcal{S}_e(\bx, t) = 0 \) for all \( t > T_s \). The forward operator \( \mathcal{A}_r \) maps the source to the resulting pressure field measured on the boundary \( \partial \Omega \) over a finite time, i.e., \( y_r \in C_0^{\infty} \left(  \partial \Omega \times [0, T] \right) \). (Note that the actual support of \( y_r \) is confined to the disk-shaped subsurface represented by \( \partial \nu_r \cap \partial \Omega \).)

\subsection{Emission Operator}
This section begins by defining Neumann-type and Dirichlet-type emitted data supported on disk-shaped subsurfaces of \( \partial \Omega \). We then introduce operators that map these surface-supported data to regularized, near-surface volumetric sources appearing on the right-hand sides of the monopole and dipole wave equations, given in Eqs.~\eqref{eq:wave-monopole} and~\eqref{eq:wave-dipole}, respectively. To that end, we start by formally defining the surface-supported source terms.
(Throughout, the subscript \( e \) denotes the fields generated by an excitation applied over the surface \( \partial \nu_e \cap \partial \Omega \).) As mentioned above, the emission process has been excluded from the definition of our forward operator (to be used within the inverse-problem framework). However, its connection with the reception process is particularly important, as it forms the basis for defining the adjoint of the reception operator, which in turn plays a key role in solving the associated minimization problem.

\begin{definition}
The \emph{Neumann-type emitted data} associated with emitter \(e\) is denoted by 
\(z_e^N \in C_0^\infty(\partial \Omega \times [0,T])\), and is defined as
\begin{align}  \label{eq:neumann-source1}
z_e^N(\bx_s, t) := 
\begin{cases}
\displaystyle \frac{\partial p_e}{\partial \bn_{e,-}}(\bx_s, t), 
& \text{if } (\bx_s, t) \in (\partial \nu_e \cap \partial \Omega) \times [0, T_s], \\[1ex]
0, & \text{otherwise},
\end{cases}
\end{align}
where \(\tfrac{\partial p_e}{\partial \bn_{e,-}}(\bx_s, t)\) denotes the inward 
normal derivative of the pressure field \(p_e\), restricted to the surface of 
emitter \(e\), with \(\bn_{e,-}\) being the inward-pointing unit normal vector 
to \(\partial \nu_e \cap \partial \Omega\). (Here, we have also used \(T_s \leq T\).)
\end{definition}
\hfill \(\square\)

\begin{definition}
The \emph{Dirichlet-type emitted data} associated with emitter \( e \) is denoted by 
\( z_e^D \in C_0^\infty(\partial \Omega \times [0,T]) \), and is defined as
\begin{align} \label{eq:Dirichlet-source1}
z_e^D(\bx_s, t) := 
\begin{cases}
p_e(\bx_s, t), 
& \text{if } (\bx_s, t) \in (\partial \nu_e \cap \partial \Omega) \times [0, T_s], \\[1ex]
0, & \text{otherwise},
\end{cases}
\end{align}
where \( p_e(\bx_s, t) \) denotes the pressure field restricted to the surface of emitter \( e \).
\end{definition}
\hfill \(\square\)

\begin{definition} \label{def:z}
We define an \emph{extension operator} that, at each arbitrary but fixed time \( t \in [0, T] \), maps emitted data \( z_e \), compactly supported on \( \partial \Omega \), to an equivalent near-surface field contained in free space. For each emitter \( e \), this operator is defined as
\begin{align}  \label{eq:z}
        \mathcal{Z}_e & : C_0^{\infty} \big( \partial \Omega \times [0, T] \big) \rightarrow C_0^\infty \big( \mathbb{R}^d \times [0, T] \big), \\
        \mathcal{Z}_e [ z_e ] (\bx, t) & := 
            -2 \displaystyle\int_{\partial \nu_e \cap \partial \Omega} dS(\bx_s)\, \delta_b(\bx - \bx_s)\, z_e(\bx_s, t), 
\end{align}
where \( \delta_b \) denotes a regularized delta distribution, whose effective support lies within a narrow region around \( \partial \Omega \), and decays smoothly away from it, as defined in Definition~\ref{def:1}.
\end{definition}
\hfill\(\square\)

From the regularized wave equation \eqref{eq:wave3}, with a regularized source term defined by the right-hand side of Eq. \eqref{eq:wave-monopole}, at an arbitrary but fixed time \(t \in [0,T]\), the near-surface field \(\mathcal{Z}_e[z_e^N]\) is added to the div-grad field \(\rho_0 \nabla \cdot \left( \frac{1}{\rho_0} \nabla p_e \right)\). 
It follows that in the limit \(b \rightarrow 0^+\), for each fixed time \(t \in [0,T]\), there is a singularity in the field 
\(\rho_0 \nabla \cdot \left( \frac{1}{\rho_0} \nabla p_e \right)\) across \(\partial \nu_e \cap \partial \Omega\). However, for positive \(b\), when a regularized Dirac delta distribution—constructed 
as a tensor product of sinc functions in Cartesian coordinates—is used, the jump \(2 \frac{\partial p_e}{\partial \bn_{e,-}}\) across the boundary does not induce a singularity. Instead, the singularity is smeared over a narrow region surrounding the boundary, with its magnitude decaying smoothly with increasing distance from 
\(\partial \nu_e \cap \partial \Omega\).

Additionally, from the equation of motion in the regularized wave system \eqref{eq:stepping3} and the force source definition \eqref{eq:force}, at an arbitrary but fixed time $t \in [0,T]$, the near-surface field \(\mathcal{Z}_e \left[ z_e^D  \right] \) is added to the normal-derivative field \( \frac{\partial p_e}{\partial \bn_{e,-}} \). (Note that we have assumed that \( \nu_r^+ \) contains the interior region \( \Omega^- \).) It follows that in the limit as $b \rightarrow 0^+$, for each fixed time \( t \in [0,T] \), there is a singularity in the field \( \frac{\partial p_e}{\partial \bn_{e,-}} \) across \(  \partial \nu_e \cap \partial \Omega   \). However, for positive $b$ when a regularized Dirac delta distribution---such as one defined in Eq.~\eqref{eq:delta-regularized}---is used, the jump \(2 p_e\) across the boundary does not induce a singularity. Instead, it is smeared over a narrow region surrounding the boundary, with its magnitude decaying smoothly with increasing distance from \(  \partial \nu_e \cap \partial \Omega  \).

\begin{definition}
We define the \emph{free-space operators} \( \mathcal{G}_e^N \) and \( \mathcal{G}_e^D \), acting on a field in free space as:
\begin{align}
    \begin{split}
        \mathcal{G}_e^N : C_0^{\infty}\left( \mathbb{R}^d \times [0,T] \right) &\rightarrow C_0^{\infty}\left( \mathbb{R}^d \times \mathbb{R}^+ \right), \quad \mathcal{G}_e^N[f](\bx, t) := f(\bx, t), \\
        \mathcal{G}_e^D : C_0^{\infty}\left( \mathbb{R}^d \times [0,T] \right) &\rightarrow C_0^{\infty}\left( \mathbb{R}^d \times \mathbb{R}^+ \right), \quad \mathcal{G}_e^D[f](\bx, t) := \frac{\partial f}{\partial \bn_{e,-}}(\bx, t).
    \end{split}
\end{align}
Here, $\mathcal{G}_e^N$ denotes the Neumann-type free-space operator associated with emitter $e$, which extends the field to the full time by setting it to zero for all $t > T$. Likewise, $\mathcal{G}_e^D$ denotes the Dirichlet-type free-space operator for emitter $e$, which first applies the inward normal derivative with respect to the infinite plane $\partial \nu_e$ and subsequently extends the resulting field to all times.
\end{definition}
\hfill \(\square\)

Having defined the surface-supported emitted data and the extension operator, we now introduce the \emph{emission operators}. These operators yield regularized representations of the source terms appearing in the regularized monopole and dipole wave equations~\eqref{eq:wave-monopole} and~\eqref{eq:wave-dipole}. These equations are subsequently recast into the wave-equation system~\eqref{eq:stepping3}, which governs wave propagation in the full space-time domain.

\begin{definition}
We define the \emph{emission operators} as mappings from surface-supported emitted data on the boundary \( \partial \Omega \) to regularized volumetric sources \( \mathcal{S}_e \), whose effective support lies within a narrow region surrounding \( \partial \Omega \) and decay smoothly in free space as the distance from the surface increases. Formally, the emission operator associated with emitter \( e \) is given by
\begin{align}  \label{eq:emission}
    & \mathcal{Z}_{\mathcal{G},e} : C_0^{\infty} \big( \partial \Omega \times [0,T] \big) \rightarrow C_0^{\infty} \left( \mathbb{R}^d \times \mathbb{R}^+ \right), \\
    & \mathcal{Z}_{\mathcal{G},e}^{\{N,D\}} \big[ z_e^{\{N,D\}} \big] := \mathcal{G}_e^{\{N,D\}} \circ \mathcal{Z}_e \big[ z_e^{\{N,D\}} \big] = \mathcal{S}_e^{\{N,D\}},
\end{align}
where we recall that \( z_e^N \) and \( z_e^D \) denote the Neumann-type and Dirichlet-type surface-supported emitted data, \( \mathcal{Z}_e \) is the extension operator, and \( \mathcal{G}_e^N \) and \( \mathcal{G}_e^D \) are the Neumann-type and Dirichlet-type free-space operators, respectively.
\end{definition}

\subsection{Reception Operator}

Leveraging the emission operator defined in the previous section, we now introduce the \emph{reception operator}, which maps the pressure wavefield in free space-time domain back to the measured data supported on the boundary \( \partial \Omega \) and within the finite time interval \( [0, T] \). To that end, we begin by formally defining the surface-supported measured data.

\begin{definition} \label{def:measured-data-Neumann}
The \emph{Neumann-type measured data} associated with receiver \( r \) and emitter \( e \) is denoted by \( y_{r,e} ^N \in C_0^{\infty} \big( \partial \Omega \times [0,T] \big) \), and is defined as
    \begin{align}
        y_{r,e}^N(\bx_s, t) := 
        \begin{cases}
            \displaystyle\frac{\partial p_e}{\partial \bn_{r,+}}(\bx_s, t), & \text{if } (\bx_s, t) \in \big( \partial \nu_r \cap \partial \Omega \big) \times [0,T], \\[0.5em]
            0, & \text{otherwise},
       \end{cases}
    \end{align}
where \( \bn_{r,+} \) denotes the outward-pointing unit normal vector to the surface \( \partial \nu_r \cap \partial \Omega \). Note that \( \bn_{r,+} \) is directed toward the exterior half-space \( \nu_r^- \).
\end{definition}
\hfill\(\square\)

\begin{definition}  \label{def:measured-data-Dirichlet}
    The \emph{Dirichlet-type measured data} associated with receiver \( r \) and emitter \( e \) is denoted by \( y_{r,e}^D \in C_0^{\infty} \big( \partial \Omega \times [0,T] \big) \), and is defined as
    \begin{align}
        y_{r,e}^D(\bx_s, t) := 
        \begin{cases}
            \displaystyle p_e(\bx_s, t), & \text{if } (\bx_s, t) \in \big( \partial \nu_r \cap \partial \Omega \big) \times [0,T], \\[0.5em]
            0, & \text{otherwise}.
        \end{cases}
    \end{align}
\end{definition}
\hfill\(\square\)

\begin{definition}    \label{def:r}
We define a \emph{restriction operator} that, at each arbitrary but fixed time \( t \in [0, T] \), maps a field \( f \) in free space to the measured data \( y \) supported on \( \partial \Omega \). For each receiver \( r \), this operator is given by
\begin{align} \label{eq:r}
    \begin{split}
        &\mathcal{R}_r : C_0^{\infty} \big( \mathbb{R}^d \times [0,T] \big) \to C_0^\infty \big( \partial \Omega \times [0, T] \big), \\
        &\quad \mathcal{R}_r \left[ f \right](\bx_s, t) = y_r(\bx_s, t) :=
        \begin{cases}
            \displaystyle \frac{1}{2}
            \int_{\mathbb{R}^d} d\bx \ \delta_{b, \partial \nu_r}(\bx_s - \bx) \, f(\bx,t), & \text{if } (\bx_s, t) \in (\partial \nu_r \cap \partial \Omega) \times [0, T], \\
            0, & \text{otherwise}.
        \end{cases}
    \end{split}
\end{align}

\noindent
Here, \( \delta_{b,\partial \nu_r}(\bx_s - \bx) \) denotes a regularized Dirac delta distribution supported on the \((d{-}1)\)-dimensional infinite hyperplane \( \partial \nu_r  \). It is defined analogously to the standard Dirac delta distribution \( \delta_b \) in \( \mathbb{R}^d \), but with support restricted to the infinite plane (units: \( \mathrm{m}^{1-d} \)). It satisfies the sifting property
\begin{align} \label{eq:sifting-surface}
    f(\bx_s) = \lim_{b \to 0^+} \int_{\partial \nu_r } dS(\bx') \, \delta_{b,\partial \nu_r}(\bx_s - \bx') \, f(\bx'), \quad \bx_s \in \partial \nu_r \cap \partial \Omega.
\end{align}
The explicit form of the surface-restricted Dirac delta distribution is derived in the following Lemma.
\end{definition}
\hfill\(\square\)

\begin{lemma} \label{lemma:extraction}
On the receiver surface, in the limit $b \to 0^+$, the field $f \in C_0^{\infty}\big(\mathbb{R}^d \times [0,T]\big)$ can be extracted by applying the restriction operator to its outward normal derivative in a narrow region surrounding the receiver surface. Specifically,
\begin{align}
    f(\bx_s, t)
    = \mathcal{R}_r \Big[ \frac{\partial f}{\partial \bn_{r,+}} (\bx, t) \Big],
    \qquad
    (\bx_s, t) \in \big(\partial \nu_r \cap \partial \Omega\big) \times [0,T],
\end{align}
where $\bn_{r,+}$ is the outward-pointing unit normal vector on $\partial \nu_r \cap \partial \Omega$, and $\mathcal{R}_r$ is defined in~\eqref{eq:r}. Moreover,
\[
    \delta_{b,\partial \nu_r}(\bx_s - \bx')
    = a_R(\bx_s) \, \delta_b(\bx_s - \bx'),
\]
where \( \delta_{b,\partial \nu_r} \) denotes the surface-restricted mollifier of bandwidth \( b \), and \( a_R(\bx_s) = \frac{d \bx_s}{dS(\bx_s)} \) is a scaling factor, which simplifies to \( a_R = b \) if one of the Cartesian coordinates coincides with the normal direction to the receiver surface.

\end{lemma}

\begin{proof}
At an arbitrary but fixed time \( t \in [0,T] \), the functions \( h, f \in C_0^{\infty} \big( \mathbb{R}^d \times [0,T] \big) \) satisfy:
\begin{align}  \label{eq:int-reception-n}
\begin{split}
     \int_{\partial \nu_r} dS \, h  f = \int_{ \nu_r^\pm} d\bx \, h \, \frac{\partial f}{\partial \bn_{r,\pm}} + \int_{\nu_r^\pm} d\bx \, \frac{\partial h}{\partial \bn_{r,\pm}}  f,
\end{split}
\end{align}
where $\bn_{r,+}$ and $\bn_{r,-}$ denote the outward and inward unit normal vectors, respectively, on the surface $\partial \nu_r \cap \partial \Omega$. We then choose the test function as
\begin{align}  \label{eq:delta-surface}
h(\bx) = \delta_b(\bx_s - \bx) +  \delta_b(\bx_s - \btx),
\end{align}
where \(\btx\) denotes the image point of \(\bx\) with respect to the infinite plane.

\noindent
Substituting Eq.~\eqref{eq:delta-surface} into both sides of 
Eq.~\eqref{eq:int-reception-n}, and summing the integral contributions from the two half-spaces, the second integral on the right-hand side vanishes because \(\tfrac{\partial h}{\partial \bn_{r,\pm}}=0 \).

Additionally, in the specific setting that one Cartesian coordinate coincides with normal direction to \( \partial \nu_r \), and by 
decomposing the regularized Dirac delta distribution into sinc functions in the coordinates tangent and normal to the infinite plane in the left-hand side integral expression, and applying its 
sifting property in the limit \(b \to 0^+\) to localize the field and its 
normal derivative from the left-hand side integral expressions, we obtain
\begin{align} \label{eq:surface-integral2}
  \frac{1}{b} \,
  \Big[ f \big|_{\partial \nu_r^+} (\bx_s, t)- f \big|_{\partial \nu_r^-} (\bx_s, t)\Big]
  \;=\;\int_{\mathbb{R}^d} d\bx \ \delta_b(\bx_s - \bx) \, 
  \frac{\partial f}{\partial \bn_{r,+}}(\bx, t), \quad (\bx_s, t) \in \big( \partial \nu_r \cap \partial \Omega \big) \times [0,T],
\end{align}
which, in the limit \( b \to 0^+ \), indicates the emergence of a 
distributional singularity in 
\( \tfrac{\partial f}{\partial \bn_{r,+}} \) 
due to the jump of \( 2 f \) across the boundary, leading to the symmetric relation
\begin{align} \label{eq:surface-integral0}
  \lim_{b \to 0^+} \, \frac{2}{b}\, f\big|_{\partial \nu_r} (\bx_s, t)
  \;=\;\int_{\mathbb{R}^d} d\bx \ \delta_b(\bx_s - \bx) \, 
  \frac{\partial f}{\partial \bn_{r,+}}(\bx, t), \quad (\bx_s, t) \in \big( \partial \nu_r \cap \partial \Omega \big) \times [0,T].
\end{align}
For a general Cartesian coordinate system, performing the change of variables \( dS(\bx_s) \to d \bx_s \) in the integral on the left-hand side of Eq.~\eqref{eq:int-reception-n} yields the same formula as in Eq.~\eqref{eq:surface-integral0}, but with \( b \) replaced by the generalized scaling factor \( a_R(\bx_s) = \frac{d \bx_s}{dS(\bx_s)} \). This yields
\begin{align} \label{eq:surface-integral}
  \lim_{b \to 0^+} \, \frac{2}{a_R(\bx_s)}\, f\big|_{\partial \nu_r} (\bx_s, t)
  \;=\;\int_{\mathbb{R}^d} d\bx \ \delta_b(\bx_s - \bx) \, 
  \frac{\partial f}{\partial \bn_{r,+}}(\bx, t), \quad (\bx_s, t) \in \big( \partial \nu_r \cap \partial \Omega \big) \times [0,T].
\end{align}

\noindent
Eq.~\eqref{eq:surface-integral} yields the action of the restriction operator 
defined in Definition~\ref{def:r} on \(\tfrac{\partial f}{\partial \bn_{r,+}}\), 
where, in Eq.~\eqref{eq:r}, the surface-restricted Dirac delta distribution is 
given by
\begin{align}  \label{eq:dirac-delta-surface}
    \delta_{b,\partial \nu_r}(\bx_s - \bx) 
    = a_R(\bx_s)  \, \delta_b(\bx_s - \bx),\quad \bx_s \in \partial \nu_r \cap \partial \Omega ,
\end{align}
with \( a_R \) serving as a scaling factor representing the relative integral measure \( \frac{d \bx_s}{dS(\bx_s)} \). In the derivation above, we assumed the receiver surface to be an infinite plane. When a quantity is restricted to an infinite plane as a smeared jump—equivalently, as a smeared singularity in the normal derivative of the pressure—the restriction operator smooths the jump on the surface, i.e., it yields non-unique values that depend on the choice of \( a_R \). This non-uniqueness is related to the singularity that arises when inverting the double-layer potential integral formula in boundary element methods. From a numerical point of view, modeling the wave propagation and solving the forward problem uniquely yields the implicit formula~\eqref{eq:surface-integral}, rather than the explicitly surface-restricted quantity, since the latter depends on the factor \( a_R \).

\end{proof}

\noindent
Conceptually, at a boundary point, the restriction operator acts as a left inverse of the extension operator introduced in Definition~\ref{def:z}, in accordance with the principle of reciprocity, which is fundamental in physics. This will be established in the following lemma.

\begin{lemma} \label{lemma:inverse-re}
For the special case \( r = e \), the restriction operator \( \mathcal{R}_r \) is a left inverse of the extension operator \( \mathcal{Z}_e \), satisfying
\[
\mathcal{R}_r \circ \mathcal{Z}_e = \mathrm{Id}.
\] 
\end{lemma}

\begin{proof}
Let \( f \) be a test function compactly supported on \( \partial \nu_e \cap \partial \Omega \). We aim to show that for \( r = e \),
\[
f =  \lim_{b \to 0^+} \,  \mathcal{R}_r \circ \mathcal{Z}_e \left[ f \right].
\]

\noindent
At each arbitrary but fixed time \( t \in [0,T]\), by substituting the expression for the extension operator \( \mathcal{Z}_e \), introduced in Definition~\eqref{def:z}, into the restriction operator \( \mathcal{R}_r \), introduced in Definition~\eqref{def:r}, and taking the limit \( b \to 0^+ \), we obtain:
\begin{align*}
\begin{split}
   \lim_{b \to 0^+}  \,  \mathcal{R}_r \circ \mathcal{Z}_e \left[ f(\bx_s, t) \right]  
    & =   \lim_{b \to 0^+} \ \frac{1}{2} \int_{\mathbb{R}^d} d\bx \,  \delta_{b,\partial \nu_r}(\bx_s - \bx)
   \ \Big[ -2 \int_{\partial \nu_e \cap \partial \Omega} dS(\bx'_s) \, \delta_b(\bx - \bx'_s) \, f(\bx'_s, t) \Big] \\
    & =  \lim_{b \to 0^+} \frac{1}{2}  \int_{\mathbb{R}^d} d\bx \,  \delta_{b,\partial \nu_r}(\bx_s - \bx)  \frac{\partial f(\bx, t)}{\partial \bn_{e,+}}  =   f(\bx_s, t),
 \end{split}
\end{align*}
for all \((\bx_s, t) \in (\partial \nu_r \cap \partial \Omega) \times [0, T]\). 
Note that, since the unit normal vector \(\bn_e\) is outward-pointing, the minus 
sign in the extension operator is dropped. Since \( r = e \), we have \( \partial \nu_r = \partial \nu_e \) and \(\bn_r = \bn_e \), and the above expression recovers \( f \). This completes the proof.

\end{proof}

\begin{definition}   \label{def:rff}
We define the \emph{free-space operators} \( \mathcal{G}_r^N \) and \( \mathcal{G}_r^D \), acting on a field in free space as:
\begin{align}  \label{eq:fdn}
    \begin{split}
        \mathcal{G}_r^N : C_0^{\infty}\left( \mathbb{R}^d \times [0,T] \right) &\rightarrow C_0^{\infty}\left( \mathbb{R}^d \times [0,T] \right), \quad \mathcal{G}_r^N[f](\bx, t) := \rho_0 \nabla \cdot \left( \frac{1}{\rho_0} \nabla f (\bx, t) \right), \\
        \mathcal{G}_r^D : C_0^{\infty}\left( \mathbb{R}^d \times [0,T] \right) &\rightarrow C_0^{\infty}\left( \mathbb{R}^d \times [0,T] \right), \quad \mathcal{G}_r^D[f](\bx, t) := \frac{\partial f}{\partial \bn_{r,+}}(\bx, t),
    \end{split}
\end{align}
where \( \bn_{r,+} \) denotes the outward-pointing unit normal vector on the surface \( \partial \nu_r \cap \partial \Omega \).

\end{definition}
\hfill\(\square\)

\begin{definition}  \label{def:r-nd}
We define the \emph{reception operator} as a mapping from a pressure wavefield in free space to surface-supported measured data on the boundary \( \partial \Omega \). Formally, the reception operator associated with receiver \( r \) is defined by
\begin{align} 
\begin{split}
    &\mathcal{R}_{\mathcal{G},r}^{\{N,D\}} : C_0^{\infty}\left( \mathbb{R}^d \times [0,T] \right) \rightarrow C_0^{\infty}\left( \partial \Omega \times (0,T) \right), \\
    &\mathcal{R}_{\mathcal{G},r}^{\{N,D\}} [f](\bx_s, t) = \left( \mathcal{R}_r \circ \mathcal{G}_r^{\{N,D\}} \right)\left[f(\bx, t)\right].
\end{split}
\end{align}
Here, \( \mathcal{R}_{\mathcal{G},r}^{\{N,D\}} \) denotes a composite operator involving \( \mathcal{G}_r^{\{N,D\}} \)—the Neumann-type and Dirichlet-type free-space operators—and \( \mathcal{R}_r \), the restriction operator. These constituent operators were introduced in Definitions~\ref{def:rff} and~\ref{def:r}, respectively.
\end{definition}
\hfill\(\square\)

\begin{lemma}  \label{lemma:extraction2}
Let \( p_e \in C_0^{\infty}\left( \mathbb{R}^d \times [0,T] \right) \) be the pressure wavefield generated by the regularized source \( \mathcal{S}_e \in C_0^{\infty}\left( \mathbb{R}^d \times \mathbb{R}^+ \right) \), and let \( y_{r,e}^{\{N,D\}} \in C_0^{\infty} \left( \partial \Omega \times [0, T] \right) \) denote the measured Neumann-type and Dirichlet-type data at receiver \(r\) corresponding to emitter \(e\). Then the associated operators \( \mathcal{R}_{ \mathcal{G},r}^{\{N,D\}} \) satisfy
\begin{align}
   y_{r,e}^{\{N,D\}} = \mathcal{R}_{\mathcal{G},r}^{\{N,D\}} \left[ p_e \right].
\end{align}
\end{lemma}

\begin{proof}
The claim follows directly from Lemma~\ref{lemma:extraction} and Definitions~\ref{def:r}, \ref{def:rff}, and~\ref{def:r-nd}. To derive the more practical Dirichlet-type measured data \(y_{r,e}^D\), we proceed exactly as in Lemma~\ref{lemma:extraction}, by setting \(f := p_e\) and choosing \(h\) as in equation~\ref{eq:delta-surface}. 

\noindent
For the extraction of the Neumann-type measured data $y_{r,e}^N$, Lemma~\ref{lemma:extraction} and its proof require slight modification; we omit the full details for brevity. In this case, we start from Eq.~\eqref{eq:int-reception-n}, setting 
\[
f := \frac{1}{\rho_0} \frac{\partial p_e}{\partial \bn_{r,+}},
\]
and defining the test function the same as in Eq.~\eqref{eq:delta-surface}. As the normal derivative of $h$ vanishes, the second term on the right-hand side of Eq.~\eqref{eq:int-reception-n} drops out. Consequently, we obtain:
\begin{align}  \label{eq:surface-inregral0}
  \lim_{b \to 0^+} \, \frac{1}{a_R(\bx_s)}\, \Big[ \frac{1}{\rho_0} \frac{\partial p_e}{\partial \bn_{r,+}} \Big|_{\partial \nu_r^+} (\bx_s, t)  - \frac{1}{\rho_0} \frac{\partial p_e}{\partial \bn_{r,+}}\Big|_{\partial \nu_r^-} (\bx_s, t) \Big]  \;=\;   \int_{\mathbb{R}^d} d\bx \ \delta_b(\bx_s - \bx)  \nabla \cdot \Big( \frac{1}{\rho_0} \nabla p_e  ( \bx , t) \Big) ,
\end{align}
which, in the limit \( b \to 0^+ \), indicates the emergence of a distributional singularity in 
\( \rho_0 \, \nabla \cdot \!\big( \tfrac{1}{\rho_0} \nabla p_e \big) \)  due to the jump of \( 2 \tfrac{\partial  p_e}{\partial \bn_{r,+}} \) across the boundary, leading to the symmetric relation
\begin{align} \label{eq:surface-integraln}
\lim_{b \to 0^+} \, \frac{2}{a_R(\bx_s)}\,  \tfrac{\partial  p_e}{\partial \bn_{r,+}} \Big|_{\partial \nu_r}  (\bx_s, t) \;=\;  \int_{\mathbb{R}^d} d\bx \ \delta_b(\bx_s - \bx)  \rho_0 \nabla \cdot \Big( \frac{1}{\rho_0} \nabla p_e( \bx , t)  \Big).
\end{align}
This yields \( y^N \) on the boundary as the action of the restriction operator, defined in Eq.~\eqref{eq:r} of Definition~\ref{def:r}, applied to the div--grad field, where the surface-restricted Dirac delta distribution is given in Eq.~\eqref{eq:dirac-delta-surface}. We remind the reader that what is uniquely determined is the implicit formula~\eqref{eq:surface-integraln}, rather than the explicit boundary quantity \( y^N := \left. \tfrac{\partial p_e}{\partial \bn_{r,+}} \right|_{\partial \nu_r} \).

\end{proof}

\noindent
By assuming the infinite plane, it follows that, in the limit $b \rightarrow 0^+$, for each fixed time $t \in [0,T]$, the free-space field 
\(\rho_0 \nabla \cdot \left( \frac{1}{\rho_0} \nabla p_e \right)\) (resp. \(\frac{\partial p_e}{\partial \bn_{r,+}}\)) exhibits a singularity across $\partial \nu_r \cap \partial \Omega$. However, for positive finite $b$, when a regularized Dirac delta distribution is employed, the jump \(2\frac{\partial p_e}{\partial \bn_{r,+}}\) (resp. \(2 p_e \)) across the boundary does not induce a singularity. Instead, the singularity is smeared over a narrow region surrounding the boundary, with its magnitude decaying smoothly with increasing distance from $\partial \nu_r \cap \partial \Omega$. Lemmas~\ref{lemma:extraction} and~\ref{lemma:extraction2} indicate that, at each time $t \in [0, T]$, the Neumann-type data $y_r^N$ (resp. the Dirichlet-type data $y_r^D$) is extracted on the surface of receiver $r$ as this smeared singularity.

\subsection{Forward Operator}

Having defined the emission and reception processes, we now introduce the forward operator in the context of ultrasound tomography. To that end, we first define the \emph{wave propagation operator}.

\begin{definition} \label{def:wave-propagation-operator}
The \textit{wave propagation operator} \(\mathcal{P}\) maps a regularized 
source \(\mathcal{S}\) in free space--time to the resulting pressure wavefield 
in free space over the interval \([0,T]\). Specifically,
\begin{align} \label{eq:wave-map}
\mathcal{P}[\mathcal{S}]: C_0^{\infty}\!\left( \mathbb{R}^d \times \mathbb{R}^+ \right) 
\;\longrightarrow\; C_0^{\infty}\!\left( \mathbb{R}^d \times [0,T] \right),
\end{align}
where \(\mathcal{P}[\mathcal{S}]\) denotes the causal solution \(p\) of 
equation~\eqref{eq:wave3} corresponding to the regularized source \(\mathcal{S}\).
\end{definition}
\hfill \(\square\)

\begin{definition}   \label{def:operator-forward}
The \textit{forward operator} maps a free space-time regularized source \( \mathcal{S}_e \) to the corresponding Neumann-type or Dirichlet-type measured data \( y_{r,e}^{\{N,D\}} \), respectively, supported on \( \partial \nu_r \cap \partial \Omega \). This map is defined as
\begin{align} \label{eq:operator-forward}
\begin{split}
    \mathcal{A}_r^{\{N,D\}} : C_0^{\infty} \big( \mathbb{R}^d \times \mathbb{R}^+ \big) \rightarrow C_0^{\infty} \big( \partial \Omega \times [0,T] \big), \\
    \mathcal{A}_r^{\{N,D\}}[\mathcal{S}_e] = \left( \mathcal{R}_{\mathcal{G},r}^{\{N,D\}} \circ \mathcal{P} \right)[\mathcal{S}_e] = y_{r,e}^{\{N,D\}},
\end{split}
\end{align}
where \(\mathcal{P}\) is the wave propagation operator, and \(\mathcal{R}_{r,\mathcal{G}}^{\{N,D\}}\) maps the resulting pressure wavefield in free space to the Neumann-type or Dirichlet-type boundary data, respectively. The quantity \( y_{r,e}^{\{N,D\}} \) corresponds to the action of \(\mathcal{A}_r^{\{N,D\}}\) on the regularized source \(\mathcal{S}_e\), which is induced by a boundary excitation on \(\partial \nu_e \cap \partial \Omega\). Repeating this process for all emitter and receiver indices \(e\) and \(r\), respectively, yields the complete dataset.
\end{definition}
\hfill \(\square\)

\section{Adjoint and Time-reversal Operators}  \label{sec:adjoint}

Given the forward operator defined above, the inverse problem seeks to minimize the objective function
\begin{align}
\sum_r \left\| \mathcal{A}_r^{\{N,D\}}[\mathcal{S}_e] - (y_{r,e}^{\{N,D\}})_{\text{measured}} \right\|_{L^2\big( \partial \Omega \times [0,T] \big)}^2
\end{align}
for all emitter indices \(e\).

Here, \( (y_{r,e}^{\{N,D\}})_{\text{measured}} \) denotes the measured surface-restricted Neumann-type or Dirichlet-type data, respectively. The unknown parameter, which serves as the minimizer of the objective function, can be the source or any parameters in the forward operator \( \mathcal{A}_r^{\{N,D\}} \), such as the acoustic properties of the medium.

To remove the singularity of the reception operator derived in Lemmas~\ref{lemma:extraction} and~\ref{lemma:extraction2}, we use the scaled objective function
\begin{align}
\sum_r \left\| \frac{1}{a_R} \mathcal{A}_r^{\{N,D\}}[\mathcal{S}_e] - \frac{1}{a_R} \bigl(y_{r,e}^{\{N,D\}}\bigr)_{\text{measured}} \right\|_{L^2\big(\partial \Omega \times [0,T]\big)}^2
\end{align}
for all emitter indices \(e\). Here, \( \frac{1}{a_R} \mathcal{A}_r^{\{N,D\}} \) denotes the scaled variant of the forward operator, i.e., Eq.~\eqref{eq:operator-forward} in which the bandwidth parameter \(a_R\) in the surface‑restricted regularized Dirac delta distribution appearing in the reception operator is canceled by this scaling, making the boundary data independent of this bandwidth parameter. (See Lemma~\ref{lemma:extraction}.) This scaling factor acts element-wise on the boundary data on the receiver surfaces.

The minimization process employs step directions that depend on both the forward operator and its adjoint. Therefore, in conjunction with the integral formulations governing the emission and reception processes of finite-sized acoustic transducers, it is essential to derive the corresponding adjoint operator under the same assumptions. The resulting formulation follows directly from the integral representations established in the preceding sections.

\subsection{Adjoint Operator}  
\label{sec:adjoint-op}
Having defined the forward operator—comprising the wave propagation and reception processes—we now proceed to derive the corresponding adjoint operator.

\subsubsection{Adjoint of the Wave Propagation Operator} 

We begin by deriving the adjoint of the wave propagation operator introduced in Definition \ref{def:wave-propagation-operator}.

\begin{lemma}    \label{lemma:adjoint-propagation}
    The action of the adjoint of the Wave Propagation Operator $\mathcal{P}$ on any test function $f^{\mathcal{S}}\in C_0^{\infty} ( \mathbb{R}^d \times [0,T] ) $ is given by:
    \begin{align}   \label{eq:wp-adjoint}
        \begin{split}
           &\mathcal{P}^*: C_0^{\infty} \big( \mathbb{R}^d \times [0,T]  \big)  \rightarrow   C_0^{\infty} \big(\mathbb{R}^d \times \mathbb{R}^+ \big) ,\\
           & \mathcal{P}^* \left[f^{\mathcal{S}}\right](\bx', t' ) = \int_0^T dt \int_{\mathbb{R}^d} d \bx \  g(\bx_{\fd}, T-t'-t ) \  f^{\mathcal{S}}(\bx, T-t).
        \end{split}
    \end{align}
\end{lemma}

\begin{proof}
The operator \( \mathcal{P} \), and its adjoint \( \mathcal{P}^* \), with respect to the standard \( L^2 \) bilinear form in \( C_0^{\infty} \big( \mathbb{R}^d \times \mathbb{R}^+ \big) \) and \( C_0^{\infty} \big( \mathbb{R}^d \times [0,T] \big) \), must satisfy:
\begin{align}  \label{eq:adjo1}
\int_0^T dt \int_{\mathbb{R}^d} d \bx  \   f^{\mathcal{S}}(\bx,t) \  \bigg[ \int_{\mathbb{R}^+} dt' \int_{\mathbb{R}^d} d \bx' \  g_+(\bx_{\fd}, t_{\fd} ) \mathcal{S}(\bx', t') \bigg]  
=  \int_0^T dt' \int_{\mathbb{R}^d} d \bx'  \   \mathcal{S}(\bx',t')  \
 \mathcal{P}^*\left[f^{\mathcal{S}}\right](\bx', t' ) 
\end{align}
for any \( \mathcal{S} \in C_0^{\infty} \big( \mathbb{R}^d \times \mathbb{R}^+ \big) \) and test function \( f^{\mathcal{S}} \in C_0^{\infty} \big( \mathbb{R}^d \times [0,T] \big) \).

Now, we rearrange the order of integration and perform the change of variables \( t \mapsto T - t \) on the left-hand side. Using the reciprocity of the Green's function,
\begin{align}  \label{eq:reciprocity}
    g(\bx_{\fd}, t_{\fd}) = g(-\bx_{\fd}, t_{\fd}),
\end{align}
we obtain:
\begin{align} \label{eq:adj0}
\int_{\mathbb{R}^+} dt' \int_{\mathbb{R}^d} d \bx'  \   \mathcal{S}(\bx', t') \bigg[  \int_0^T dt \int_{\mathbb{R}^d} d \bx  \  g(\bx_{\fd}, T - t' - t ) \  \ f^{\mathcal{S}}(\bx, T - t)  \bigg].
\end{align}

Based on the right-hand side of Eq.~\eqref{eq:adjo1}, we identify the bracketed term in Eq.~\eqref{eq:adj0} as the adjoint operator \( \mathcal{P}^* \) acting on \( f^{\mathcal{S}} \), as defined in Eq.~\eqref{eq:wp-adjoint}.
\end{proof}

\subsubsection{Adjoint of the Reception Operator} 

Having derived the adjoint operator \( \mathcal{P}^* \), we now proceed to derive the adjoint of the reception operators \( \mathcal{R}_{\mathcal{G},r}^{\{N,D\}} \), which are, respectively, composite operators consisting of the Neumann-type or Dirichlet-type free-space operators introduced in Definition~\ref{def:rff}, and the restriction operator introduced in Definition~\ref{def:r}. 

\begin{lemma}   \label{lemma:rt-n}
The action of the adjoint of the scaled variant of the Neumann-type reception operator \( \frac{1}{a_R}\mathcal{R}_{\mathcal{G},r}^N \), introduced in Definition~\eqref{def:r-nd}, on any test function \( \frac{1}{a_R} f^N \in C_0^{\infty}\!\left(\partial \Omega \times [0,T]\right) \) is given by
\begin{align}       \label{eq:meas-adj}
    \begin{split}
    \left(  \frac{1}{a_R} \mathcal{R}_{\mathcal{G},r}^N \right)^* \!:&\;  C_0^{\infty} \big(\partial \Omega \times [0, T] \big) \rightarrow C_0^{\infty} \big(\mathbb{R}^d \times [0,T] \big),\\
     \left(  \frac{1}{a_R} \mathcal{R}_{\mathcal{G},r}^N \right)^*[\frac{1}{a_R} f^N](\bx, t) &= 
\dfrac{1}{2} \displaystyle\int_{\partial \nu_r \cap \partial \Omega} d S(\bx_s)\,  \nabla \cdot \left( \dfrac{1}{\rho_0} \nabla \Big( \rho_0 \delta_b(\bx - \bx_s) \Big) \right) f^N(\bx_s, t), 
    \end{split}
\end{align}
where we have used the scaled operator to remove the singularity of the reception operator.

\end{lemma}

\begin{proof}
Recall that the scaled operator \(   \frac{1}{a_R}\mathcal{R}_{\mathcal{G},r}^N  \) map the free-space field to the Neumann-type boundary data
\(   f^N \) scaled by \( \frac{1}{a_R} \). Here, \(a_R(\bx_s) = \frac{d\bx_s}{dS(\bx_s) } \), as defined in Eq.~\eqref{eq:dirac-delta-surface}. At any arbitrary but fixed time \( t \in [0, T] \), the operator \(   \frac{1}{a_R}\mathcal{R}_{\mathcal{G},r}^N  \) and its adjoint, with respect to the standard \( L^2 \) bilinear form in  \( C_0^{\infty}(\mathbb{R}^d \times [0,T]) \) and \( C_0^{\infty}(\partial \Omega \times [0,T]) \), must satisfy
\begin{align}\label{eq:meas-adj-test1}
\begin{split}
 & \int_{\partial \nu_r \cap \partial \Omega} d\bx_s \,\, \frac{1}{a_R(\bx_s)}  f^N(\bx_s,t) \frac{1}{2} \int_{\mathbb{R}^d} d\bx \ \delta_b(\bx_s - \bx)\, \rho_0(\bx) \nabla \cdot \left( \frac{1}{\rho_0(\bx)} \nabla f(\bx, t) \right) \\
= & \int_{\mathbb{R}^d} d\bx \ \left[ \left( \frac{1}{a_R} \mathcal{R}_{\mathcal{G},r}^N \right)^* [\frac{1}{a_R} f^N]\right] f(\bx, t)
\end{split}
\end{align}
for any \( f \in C_0^{\infty}(\mathbb{R}^d \times [0,T]) \) and \( \frac{1}{a_R} f^N \in C_0^{\infty}(\partial \Omega \times [0, T]) \), respectively. 

Applying integration by parts twice to the left-hand side of Eq.~\eqref{eq:meas-adj-test1} (using the compact support to discard boundary terms), and interchanging the order of integration, yield
\begin{align}  \label{eq:meas-adj-test2} 
   \int_{\mathbb{R}^d} d\bx \left[ \frac{1}{2} \int_{\partial \nu_r \cap \partial \Omega} d\bx_s \, \nabla \cdot \left( \frac{1}{\rho_0 (\bx)} \nabla \Big( \rho_0(\bx) \delta_b(\bx - \bx_s) \Big) \right) \frac{1}{a_R(\bx_s)} f^N(\bx_s, t) \right] f(\bx, t).
\end{align}
By comparing the integral expression in~\eqref{eq:meas-adj-test2} with the right-hand side of Eq.~\eqref{eq:meas-adj-test1}, we identify the bracketed term as \( \left( \frac{1}{a_R} \mathcal{R}_{\mathcal{G},r}^N \right)^* [\frac{1}{a_R} f^N] \), thereby confirming the adjoint expression in Eq.~\eqref{eq:meas-adj}.
\end{proof}

\begin{lemma}   \label{lemma:rt-d}
The action of the adjoint of the scaled variant of the Dirichlet-type reception operator \( \frac{1}{a_R} \mathcal{R}_{\mathcal{G},r}^D \), as introduced in Definition~\eqref{def:r-nd}, on any test function \(  \frac{1}{a_R} f^D \in  C_0^{\infty} \left(\partial \Omega \times [0, T] \right) \) is given by:
\begin{align}       \label{eq:meas-adj-d}
    \begin{split}
    \left(  \frac{1}{a_R} \mathcal{R}_{\mathcal{G},r}^D \right)^* \!:&\;  C_0^{\infty} \big(\partial \Omega \times [0, T] \big) \rightarrow C_0^{\infty} \big(\mathbb{R}^d \times [0,T] \big),\\
     \left(  \frac{1}{a_R} \mathcal{R}_{\mathcal{G},r}^D \right)^* [  \frac{1}{a_R} f^D ](\bx, t) &= 
     \dfrac{1}{2}   \displaystyle\int_{\partial \nu_r \cap \partial \Omega} d S(\bx_s)  \,   \nabla \delta_b(\bx - \bx_s) \cdot \big[f^D(\bx_s, t) \bn_{r,-} \big],
    \end{split}
\end{align}
where we have used the scaled operator to remove the singularity of the reception operator. Also, \( \bn_{r,-} \) is the inward-pointing unit normal vector to \( \partial \nu_r \cap \partial \Omega \).

\end{lemma}

\begin{proof}
Recall that the scaled operator \(   \frac{1}{a_R}\mathcal{R}_{\mathcal{G},r}^D  \) map the free-space field to the Dirichlet-type boundary data
\(   f^D \) scaled by \( \frac{1}{a_R} \). Here, \(a_R(\bx_s) = \frac{d\bx_s}{dS(\bx_s) } \), as defined in Eq.~\eqref{eq:dirac-delta-surface}. At any arbitrary but fixed time \( t \in [0, T] \), the operator \( \frac{1}{a_R}\mathcal{R}_{\mathcal{G},r}^D  \) and its adjoint, with respect to the standard \( L^2 \) bilinear form in \( C_0^{\infty}(\mathbb{R}^d \times [0,T]) \) and \( C_0^{\infty}(\partial \Omega \times [0,T]) \), must satisfy
\begin{align}\label{eq:meas-adj-test1-d}
\begin{split}
 & \int_{\partial \nu_r \cap \partial \Omega} d \bx_s \, \, \frac{1}{a_R} f^D(\bx_s,t)  \frac{1}{2}   \int_{\mathbb{R}^d} d\bx \  \delta_b(\bx_s - \bx)\  \big[ \nabla  f(\bx, t) \cdot \bn_{r,+}\big]   \\
= &\int_{\mathbb{R}^d} d\bx \ \left[ \left(\frac{1}{a_R} \mathcal{R}_{\mathcal{G},r}^D \right)^* [\frac{1}{a_R} f^D](\bx, t) \right] f(\bx, t)
\end{split}
\end{align}
for any  \( f \in C_0^{\infty}(\mathbb{R}^d \times [0,T]) \) and \( \frac{1}{a_R} f^D \in C_0^{\infty}(\partial \Omega \times [0, T]) \), respectively.

Applying an integration by parts to the left-hand side of Eq.~\eqref{eq:meas-adj-test1-d} (using the compact support to discard boundary terms), and interchanging the order of integration, yields
\begin{align}  \label{eq:meas-adj-test2-d} 
   \int_{\mathbb{R}^d} d\bx \left[   \frac{1}{2} \displaystyle\int_{\partial \nu_r \cap \partial \Omega} d \bx_s\,  \nabla \delta_b(\bx - \bx_s) \cdot \big[\frac{1}{a_R(\bx_s)}f^D(\bx_s, t) \bn_{r,-} \big] \right] f(\bx, t).
\end{align}

By comparing the integral expression in~\eqref{eq:meas-adj-test2-d} with the right-hand side of Eq.~\eqref{eq:meas-adj-test1-d}, we identify the bracketed term as \( \left( \frac{1}{a_R} \mathcal{R}_{\mathcal{G},r}^D \right)^* [\frac{1}{a_R} f^D] \), thereby confirming the adjoint expression in Eq.~\eqref{eq:meas-adj-d}.
\end{proof}

\begin{lemma}   \label{lemma-adjoint-operator}
The action of the adjoint of the scaled forward operator \(\frac{1}{a_R} \mathcal{A}_r^{\{ N,D\}} \), introduced in Definition~\ref{def:operator-forward}, on any scaled measured boundary data \(  \frac{1}{a_R} y_{r,e}^{ \{N,D\} } \in  C_0^{\infty} \big(\partial \Omega \times [0, T] \big) \) is given by
\begin{align}  \label{eq:adjoint-op}
    \left( \frac{1}{a_R} \mathcal{A}_r^{\{ N,D\}} \right)^*  \left[   \frac{1}{a_R} y_{r,e}^{\{N,D\}}  \right](\bx', t') = p_{e,r}^*(\bx', T - t'),
\end{align}
where \( p_{e,r}^* \) is the free space-time adjoint wavefield associated with emitter \(e\) and receiver \(r\), and is the solution to the adjoint wave equation:
\begin{align}  \label{eq:wave-adjoint}
    \left[ \frac{1}{c^2} \frac{\partial^2 }{\partial t'^2} - \rho_0  \nabla_{\bx'} \cdot \left( \frac{1}{\rho_0} \nabla_{\bx'} \right) \right] p_{e,r}^*(\bx',t') = \left( \mathcal{S}^*\right)_{e,r}^{\{N,D\}} (\bx', T - t'),
\end{align}
with Cauchy conditions:
\begin{align}   \label{eq:adjoint-initial}
    p_{e,r}^*(\bx', 0) = 0, \quad \frac{\partial}{\partial t'} p_{e,r}^*(\bx', 0) = 0, \quad  \bx' \in \mathbb{R}^d.
\end{align}
\noindent
Here, the Neumann-type and Dirichlet-type regularized adjoint sources \( \left( \mathcal{S}^* \right)_{e,r}^{\{N,D\}} \) are defined as:
\begin{align}  \label{eq:source-adjoint}
    \left( \mathcal{S}^* \right)_{e,r}^{\{N,D\}} (\bx', t') = 
    \begin{cases}
         \displaystyle \frac{1}{2} \int_{\partial \nu_r \cap \partial \Omega} dS(\bx_s) \, \nabla \cdot \left( \frac{1}{\rho_0(\bx')} \nabla \Big( \rho_0(\bx') \, \delta_b(\bx' - \bx_s) \Big) \right) y_{r,e}^N(\bx_s, t'), \\[2ex]
        \displaystyle \frac{1}{2} \int_{\partial \nu_r \cap \partial \Omega} dS(\bx_s) \, \nabla \delta_b(\bx' - \bx_s) \cdot \left[ y_{r,e}^D(\bx_s, t') \, \bn_{r,-} \right],
    \end{cases}
\end{align}
respectively.\footnote{Here we derive the adjoint operator $\left( \frac{1}{a_R}\mathcal{A}_{r}^{\{N,D\}} \right)^*$, i.e., the adjoint corresponding to a single receiver indexed by $r$. The full adjoint operator $\left( \frac{1}{a_R}\mathcal{A}^{\{N,D\}} \right)^*$ is obtained by summing over all receivers $r$ in the adjoint source formula~\eqref{eq:source-adjoint}, which yields $\sum_{r} \left(\mathcal{S}^{*}\right)_{e,r}^{\{N,D\}}$.}

\end{lemma}

\begin{proof}
The adjoint operator \(  \left( \frac{1}{a_R} \mathcal{A}_r^{\{ N,D\}} \right)^* \) is composed of the adjoints of the constituent operators:
\begin{align}
    \left( \frac{1}{a_R}\mathcal{A}_r^{\{ N,D\}} \right)^* = \mathcal{P}^* \left( \frac{1}{a_R}\mathcal{R}_{\mathcal{G}, r}^{\{N,D\}} \right)^* ,
\end{align}
where \( \mathcal{P}^* \) is given in Lemma~\ref{lemma:adjoint-propagation}, and \( \left(\frac{1}{a_R} \mathcal{R}_{\mathcal{G}, r}^{\{N,D\}} \right)^* \) are derived by Lemmas~\ref{lemma:rt-n} and~\ref{lemma:rt-d}, respectively. Composing these adjoints yields the time-reversed adjoint wave equation as stated in Eqs.~\eqref{eq:adjoint-op}--\eqref{eq:wave-adjoint}, with Cauchy conditions~\eqref{eq:adjoint-initial} and a time-reversed regularized source satisfying Eq.~\eqref{eq:source-adjoint}.
\end{proof}

\noindent
Comparing the second line of the regularized adjoint source in Eq.~\eqref{eq:source-adjoint} with the right-hand side of the regularized dipole wave equation~\eqref{eq:wave-dipole} reveals that the action of the adjoint operator on the more practical Dirichlet boundary data \( y_{r,e}^D \) is equivalent to a time-reversed variant of the dipole integral formulation~\eqref{eq:wave-dipole}. In this adjoint formulation, the regularized source---which differs from the right-hand side of the dipole wave equation~\eqref{eq:wave-dipole} only by a constant factor---is also time-reversed.

\subsection{Time-Reversal Operator}
\label{sec:tr-op}

For each receiver \( r \), the wavefield within the half-space \( \nu_r^+ \)—which contains the domain \( \Omega^- \) and the support of the source \( \mathcal{S}_e \)—can be recovered as the solution to an interior-field boundary-value problem. This solution is expressed through an integral representation involving the field and its normal derivative evaluated on the surface \( \partial \nu_r \cap \partial \Omega \).

To achieve this, we adopt the same approach as in Section~\ref{sec:gr}. In Eq.~\eqref{eq:gr-th2}, we take the limits \( t_0 \to -\infty \) and \( t_1 \to \infty \), and choose the domain of integration to be the half-space \( \nu_r^+ \). For times \( t \) greater than the source turn-off time \( T_s \)—i.e., when the source vanishes—the wavefield inside \( \nu_r^+ \), which contains \( \Omega^- \) and the support of the source \( \mathcal{S}_e \), can be expressed in terms of the field and its normal derivative on the surface \( \partial \nu_r \cap \partial \Omega \). This leads to the following \textit{time-reversed interior-field} integral equation~\cite{Devaney}:
\begin{align} \label{eq:kh-interior}
\begin{split}
   & \int_{-\infty}^{\infty} dt \int_{\partial \nu_r \cap \partial \Omega} dS \, 
    \Bigg[ 
        g_+(\bx_{\fd}, T - t' - t) \left( - \frac{\partial p_e}{\partial \bn_{r,-}}(\bx_s, T - t) \right)
        + \frac{\partial}{\partial \bn_{r,-}} g_+(\bx_{\fd}, T - t' - t) \, p_e(\bx_s, T - t)  
    \Bigg]\\
   & =  p_{e,r}^T(\bx', t'), \quad \bx' \in \nu_r^+, 
\end{split}
\end{align}
which holds for times \( t' \in (T_s, +\infty) \). Here, the change of variables \( t \rightarrow T - t \) has been applied, and we have also invoked the reciprocity of the Green’s function, as defined in Eq.~\eqref{eq:reciprocity}. Moreover, \( \bn_{r,-} \) denotes the inward-pointing unit normal to \( \partial \nu_r \cap \partial \Omega \). This boundary integral formula yields \( p_{e,r}^T \), the time-reversed wavefield corresponding to emitter \( e \) and receiver \( r \), in the half space \(\nu_r^+\).

Taking the domain of integration in Eq.~\eqref{eq:gr-th2} to be the exterior half-space \( \nu_r^- \) yields the \textit{second Helmholtz identity}, which establishes a relationship between the field and its normal derivative over \( \partial \nu_r \cap \partial \Omega \). This dependency enables the surface integral formula~\eqref{eq:kh-interior}, corresponding to the half-space \( \nu_r^+ \) and known as the \textit{first Helmholtz identity}, to be evaluated solely in terms of either the field or its normal derivative, thereby allowing the imposition of Dirichlet or Neumann boundary conditions, respectively.

Since each receiver surface is flat, the boundary data is assumed to be supported on a collection of infinite planes. This allows the Neumann and Dirichlet Green's functions to be replaced by the free-space Green's function, multiplied by a factor of 2, via the \emph{method of images} explained in Section~\ref{sec:ray-som}.

Consequently, the action of the time-reversal operator $\mathcal{A}_r^T$ on any Neumann-type or Dirichlet-type data $ y_{r,e}^{\{N,D\} }  \in C_0^{\infty} \big( \partial \Omega \times [0,T] \big) $ is given by:
\begin{align}  \label{eq:tr-operator-reversed}
    A_r^T[ y_{r,e}^{ \{N,D\} }](\bx', t') = p_{e,r}^T(\bx', t'),
\end{align}
where $p_{e,r}^T$ satisfies for any \(\bx' \in \nu_r^+\) and \(t'>T_s\)
\begin{align} \label{eq:tr-operator}
p_{e,r}^T(\bx', t')  = 
    \begin{cases}
    \begin{split}
    &2\int_0^T dt \int_{\partial \nu_r \cap \partial \Omega} dS \   g_+(\bx' -\bx_s, T-t'-t)  y_{r,e}^{N}(\bx_s, T-t),\\
      &2\int_0^T dt \int_{\partial \nu_r \cap \partial \Omega} dS \ \Big[ \frac{\partial }{\partial \bn_-} g_+(\bx'-\bx_s, T-t'-t) \Big] y_{r,e}^D (\bx_s, T-t).
    \end{split}
    \end{cases}
\end{align}
Here, \( \mathcal{A}_r^T \) denotes the \textit{time-reversal} operator associated with receiver \(r\), which acts on \( y_{r,e}^{\{N,D\}} \) and produces \( p_{e,r}^T \), the time-reversed interior wavefield corresponding to emitter \(e\) and receiver \(r\). The time-reversal operator, as defined in Eq.~\eqref{eq:tr-operator-reversed}, together with the first and second lines of Eq.~\eqref{eq:tr-operator}, represents the time-reversed, interior-field counterparts of the monopole and dipole integral formulations in Eqs.~\eqref{eq:gr-monopole} and~\eqref{eq:gr-dipole1}, respectively. Recall that a regularized version of this operator coincides with the \emph{adjoint} operator defined by 
Eqs.~\eqref{eq:adjoint-op} and~\eqref{eq:wave-adjoint}, together with the Cauchy 
initial conditions in Eq.~\eqref{eq:adjoint-initial} and the time-reversed 
regularized source in Eq.~\eqref{eq:source-adjoint}, differing only by a 
constant factor.

\section{Full Discretization of the Wave Equation on a Regular Grid}
\label{sec:wave-discretised}

Having defined the forward and adjoint operators arising from the regularized wave equation \eqref{eq:wave3}, with finite-sized emitters and receivers, this section outlines the procedure for discretizing the corresponding first-order wave equation system~\eqref{eq:stepping3} on a regular grid, with particular emphasis on the incorporation of sources in the wave equation.

Recall that the \emph{monopole} integral formula~\eqref{eq:gr-monopole} represents the pressure field in terms of an integral formula involving a monopole source, given by either \( -\frac{\partial p}{\partial \bn} \) or its equivalent \( \rho_0 \frac{\partial u^{\bn}}{\partial t} \), confined to a surface.

Similarly, the \emph{dipole} integral formula~\eqref{eq:gr-dipole1} represents the pressure field in terms of an integral formula involving a dipole source, given by  \( p \bn \), which is also confined to a surface.

\subsection{Discretized Algorithm}  
\label{sec:discretised}

Let \( \bX = \{\bX^\zeta : \zeta \in \{1, \ldots, d\} \} \) be the position of sampled points on a regular grid, where \( \zeta \) indexes the Cartesian coordinates, and let \( \Delta x \) denote the uniform grid spacing applied to all coordinate directions \( \zeta \), without loss of generality. Each sampled point is indexed by \( i \in \{ 1,\ldots N_i \} \). Furthermore, let \( \bt \in \{0, \ldots, N_t\} \) represent the discrete time steps sampled within the measurement period \( t \in \{0, \ldots, T\} \), where \( N_s \) denotes the index of sampled time step corresponding to the turn-off time of source radiation, \( T_s \).

A bar notation is used to denote fields in the fully discretized domain. The discretization of the wave equation system \eqref{eq:stepping3} on a grid staggered in both space and time is outlined in Algorithm \ref{alg:1}.

\begin{algorithm}
    \caption{Full-discretization at time step \( \bt \in \{0, \ldots, N_t-1\} \)}
    \label{alg:1}
    \begin{algorithmic}[1]
        \State \textbf{Input:} \( \bar{c}, \bar{\rho}_0, \Delta t, \Lambda^\zeta, \bar{\mathcal{S}}_f^{\zeta}(\bX, \bt) \ (\zeta \in \{1, \ldots, d\}), \bar{\mathcal{S}}_m(\bX, \bt+\frac{1}{2})  \)
        \State \textbf{Initialize:} \( \bar{p}(\bX, 0) = 0, \ \bar{\bro}(\bX, 0) = 0, \ \bar{\bu}(\bX, -\frac{1}{2}) = 0 \) \Comment{Set Cauchy conditions} 
        \For{\( \bt = 0, \ldots, N_t-1 \)}
            \State \( \bar{u}^\zeta (\bX,\bt + \frac{1}{2})  \gets \Lambda^\zeta \ \Big[\Lambda^\zeta \bar{u}^\zeta(\bX, \bt -\frac{1}{2}) - \Delta t \frac{1}{\bar{\rho}_0(\bX)} \frac{\partial}{\partial \zeta} \bar{p}(\bX, \bt)\Big] 
            + \Delta t \bar{\mathcal{S}}_f^{\zeta} (\bX, \bt) \) \Comment{Update \( \bu \)}
            \State \( \bar{\rho}^\zeta(\bX, \bt+1)  \gets \Lambda^\zeta \Big[\Lambda^\zeta \bar{\rho}^\zeta(\bX, \bt) - \Delta t \bar{\rho}_0(\bX) \frac{\partial}{\partial \zeta} \bar{u}^\zeta(\bX, \bt+\frac{1}{2})\Big] + \Delta t \bar{\mathcal{S}}_m^\zeta (\bX, \bt+\frac{1}{2}) \) \Comment{Update \( \rho \)}
            \State \( \bar{p}(\bX, \bt+1) \gets \bar{c}(\bX)^2 \sum_{\zeta = 1}^d \bar{\rho}^\zeta (\bX, \bt+1) \) \Comment{Update \( p \)}
            \State \textbf{Record} \( \bar{p}(\bX, \bt+1 ) \)
        \EndFor
        \State \textbf{Output:} \( \bar{p}(\bX, \bt) \) for \( \bt \in \{1, \ldots, N_t\} \) \Comment{Recorded pressure over \( t \in (0, T) \)}
    \end{algorithmic}
\end{algorithm}

In Algorithm~\ref{alg:1}, \( \Lambda^\zeta = e^{-\alpha^\zeta \Delta t / 2} \) denotes 
the direction-dependent Perfectly Matched Layer (PML) operator, where \( \alpha^\zeta \) 
is the virtual absorption coefficient of the PML layer along the Cartesian coordinate \( \zeta \) \cite{Mast,Tabei,Treeby}. To accommodate this directional dependency, the scalar mass source is decomposed along the Cartesian axes as 
\( \bar{\mathcal{S}}_m^\zeta = (1/d)\, \bar{\mathcal{S}}_m \), which remains identical in 
all Cartesian directions \( \zeta \) under the assumption of isotropy. Here, 
\( \bar{\mathcal{S}}_m \) is obtained by discretizing either 
Eq.~\eqref{eq:mass-u} or Eq.~\eqref{eq:mass-p2}, depending on whether a monopole or a 
monopole-like source is being modeled, respectively.

Note that, in contrast to the velocity vector \( \bar{\bu} \), the acoustic density is a scalar field; however, it is virtually treated as a vector-valued 
field \( \bar{\bro} \) to ensure compatibility with the direction-dependent PML formulation 
\cite{Mast,Tabei,Treeby}. As discussed earlier, the assumptions underlying Eq.~\eqref{eq:mass-p2} are not strictly valid; nevertheless, owing to its widespread use \cite{Treeby}, it will be employed as a benchmark to demonstrate the effectiveness of the force-based source formulation introduced in Eq.~\eqref{eq:force}.

In addition, \( \bar{\mathcal{S}}_f^\zeta \) denotes the discretized form of the \( \zeta \)-component of the vector-valued force source defined in Eq.~\eqref{eq:force}.

\subsection{Discretization of the Directional Gradients} \label{sec:grad-dis}

A k-space pseudo-spectral method is employed for the discretization of directional gradients of fields in Algorithm~\ref{alg:1} \cite{Mast,Tabei,Treeby}. It is important to note that incorporating sources in Algorithm~\ref{alg:1} remains independent of the specific method chosen for discretizing the directional gradients of fields.

\subsection{Triangulation of Emitter Surface} 
\label{sec:source-approx}

This section describes the discretization of a source supported on a finite-sized flat surface. As outlined in Section~\ref{sec:forward-def}, the emitted boundary data is supported on disjoint, finite-sized subsurfaces contained within the boundary \( \partial \Omega \). Each emitting subsurface, indexed by \( e \), is formed by the intersection of an infinite plane associated with emitter \( e \) and the boundary \( \partial \Omega \), and is denoted by \( \partial \nu_e \cap \partial \Omega \). 

Accordingly, the \((d{-}1)\)-dimensional disk-shaped subsurface \( \partial \nu_e \cap \partial \Omega \) is partitioned into a union of non-overlapping primitive-shaped elements. (For \( d=3 \), the elements are chosen as triangles. Similarly, in the case \( d=2 \), the elements reduce to lines.) These elements form a triangulation \( \mathcal{T}_e \) of the emitter surface. Each triangle \( K \in \mathcal{T}_e \) is associated with a local set of node indices \( l(K) \subseteq \{1, \dots, N_j\} \), where \( N_j \) is the total number of nodes on the emitter surface. The set \( l(K) \) identifies the \( N_{l(K)} \) nodes that define the vertices of element \( K \) (\(N_{l(K)} =d\)), and these nodes are used for interpolation and quadrature on \( K \). 

\noindent
In the sequel, we neglect the dependency on \( e \) for brevity.
 We assume that the source field \( f \) is supported on the emitter surface, and varies at the nodes \( j \) according to the relation
\[
f(\bx_j) = N_j \, a_{0,j} \, f_{\mathrm{emitter}},
\]
where \( \bx_j \) denotes the position of node \( j \), and \( a_{0,j} \) is a sensitivity coefficient satisfying
\[
\sum_{j=1}^{N_j} a_{0,j} = 1.
\]
Here, \( f_{\mathrm{emitter}} \) is a scalar value representing the strength of the emitter.

\subsection{Discretized Mass and Force Sources} \label{sec:dis-mass&force}

This section outlines the discretization of the mass and force sources, following the triangulation procedure introduced in Section~\ref{sec:source-approx}. For a source supported on a surface embedded in \( \mathbb{R}^d \), the size of each element \( K \) (i.e., the area for \( d=3 \)) is denoted by \( s_K \).

\noindent
Correspondingly, the full discretization of the surface-supported monopole source defined by Eq.~\eqref{eq:mass-u} and the monopole-like source defined by Eq.~\eqref{eq:mass-p2} at grid point \( i \) yields, respectively:
\begin{align}  \label{eq:mass-u-dis}
    \bar{\mathcal{S}}_{(m,u^{\bn})}(\bX_i, \bt) \approx  2 \sum_{K=1}^{N_K}  \frac{s_K}{N_{l(K)}}
    \sum_{\substack{j \mid  j \in l(K), \bar{\delta}_b(\bX_i- \bx_j) >  \eps }} \bar{\delta}_b(\bX_i - \bx_j)  
    \left[ \rho_0(\bx_j) \bu(\bx_j, \bt) \cdot \bn \right],
\end{align}
and
\begin{align}  \label{eq:mass-p2-dis}
    \bar{\mathcal{S}}_{(m,p)}(\bX_i, \bt) \approx 2 \sum_{K=1}^{N_K}  \frac{s_K}{N_{l(K)}}  
    \sum_{\substack{j \mid j \in l(K), \bar{\delta}_b(\bX_i- \bx_j) > \eps }} \bar{\delta}_b(\bX_i - \bx_j) 
    \left[ \frac{1}{c(\bx_j)} p(\bx_j, \bt) \right],
\end{align}
where the bandwidth parameter \(b\) in the regularized Dirac delta distribution, defined in Eq.~\eqref{eq:delta-regularized}, must be set equal to the grid spacing \(\Delta x\) to accommodate \(\delta_b\) with the grid spacing. Also, $\varepsilon$ is a scalar value chosen in the range \([0, \frac{0.05}{b^d}]\), balancing accuracy and computational cost. 

\noindent
Recall from Section~\ref{sec:doublet} that using Eq.~\eqref{eq:mass-p2-dis} to model a monopole-like surface source in terms of pressure imposes two assumptions—far-field propagation and omnidirectionality—which may not be valid in general.

\noindent
\noindent
Furthermore, the full discretization of the force source formula~\eqref{eq:force} at grid point \( i \) yields:
\begin{align}  \label{eq:force-dis}
   \bar{\boldsymbol{\mathcal{S}}}_f (\bX_i, \bt) \approx \frac{2}{\rho_0(\bX_i)}  \sum_{K=1}^{N_K} \frac{s_K }{N_{l(K)}}
   \sum_{\substack{ j \mid j \in l(K), \bar{\delta}_b(\bX_i- \bx_j) > \eps}} \bar{\delta}_b(\bX_i - \bx_j)  
   \left[p(\bx_j, \bt) \bn \right].
\end{align}

Recall from Section \ref{sec:force-source} that Algorithm~\ref{alg:1}, with a discretized force source $\bar{\boldsymbol{\mathcal{S}}}_f$ defined by Eq.~\eqref{eq:force-dis}, approximates the original dipole integral formula~\eqref{eq:gr-dipole1} without imposing any limiting assumptions.

\section{Numerical Results}  
\label{sec:num_results}  

\noindent  
This section evaluates the accuracy of the full-waveform approximation to the wave equation, as implemented in Algorithm~\ref{alg:1}. The analysis compares the approximated wavefields with analytic solutions in three representative cases:

\begin{enumerate}
  \item Modeling the action of the Green's function on a point source in the primary Green’s formula~\eqref{eq:gr-prime} using a discretized, regularized mass source localized at a single point and added to the equation of continuity in Algorithm~\ref{alg:1}.
  
  \item Approximating the monopole integral formula~\eqref{eq:gr-monopole} using a regularized mass source discretized according to Eq.~\eqref{eq:mass-u-dis}, and added to the equation of continuity in Algorithm~\ref{alg:1}.
  
  \item Modeling the dipole integral formula~\eqref{eq:gr-dipole1} using a regularized force source discretized according to Eq.~\eqref{eq:force-dis}, and added to the equation of motion in Algorithm~\ref{alg:1}.
\end{enumerate}

The evaluation begins with the action of the Green's function on a point-like volumetric radiation source and is then extended to simulate wavefields generated by finite-sized apertures, including a monopole source distributed over a disk-shaped surface and a dipole source applied over a similar region.

\subsection*{Full-waveform approach}
Wave simulations in free space were performed using the k-Wave toolbox, which implements Algorithm~\ref{alg:1} and employs a k-space pseudospectral method to compute the directional gradients of the fields~\cite{Treeby}. 
The mass and force source definitions follow the discretized formulas~\eqref{eq:mass-u-dis}--\eqref{eq:force-dis}.

A computational grid with sampled positions in the range $\big[-7.14, +7.14 \big] \times \big[-7.14, +7.14 \big] \times \big[-7.14, 0.5 \big] \text{cm}^3$ and with a spacing of 0.4~mm along all Cartesian coordinates was used. The sound speed and ambient density were set to 1540~$\text{ms}^{-1}$ and 1000~$ \text{kgm}^{-3}$, respectively, in a homogeneous medium. The maximum frequency supported by the grid for wave simulation was determined by the \textit{Shannon-Nyquist} limit \cite{Tabei}. For a homogeneous medium with sound speed $c$, the maximum supported frequency is given by \cite{Treeby}:
\begin{align} \label{eq:fmax}
    f_\text{max}= \frac{c}{2 \Delta x},
\end{align}
which equals 1.925~MHz in this experiment. The time step was set to 0.04~$ \mu \text{s}$. 

This section evaluates the full-waveform approach for approximating the action of the causal Green's function on a point source, as well as the \textit{monopole} and \textit{dipole} integral formulas derived in \eqref{eq:gr-monopole} and \eqref{eq:gr-dipole1}, respectively. First, the mass source is defined in terms of a regularized source, $\mathcal{S}$, located at a single point. Then, Eq. \eqref{eq:mass-u-dis} is used to extend the mass source formulation to model a monopole source supported on a disk-shaped surface with a radius of $8$ mm. Finally, Eq. \eqref{eq:force-dis} is used to model a vector-valued force source in terms of a dipole source, $p \bn$, also supported on a disk-shaped surface with a radius of $8$ mm. 

\subsection*{Analytical Approach}
The accuracy of the full-waveform approach was evaluated by comparing its approximated wavefields with analytical solutions, which served as benchmarks. Specifically, the action of the causal Green's function on a source confined to a single point was computed analytically using a frequency-domain variant of the primary Green's formula \eqref{eq:gr-prime}, where the 3D Green's function acts on a point source. For a source confined to a single point, the spatial integral in Eq. \eqref{eq:gr-prime} is omitted.

Furthermore, the open-source \textit{Field II} toolbox was employed to calculate the \textit{monopole} integral formula \eqref{eq:gr-monopole} and the \textit{dipole} integral formula \eqref{eq:gr-dipole1} in the time domain \cite{Jensen1,Jensen2}. These analytical solutions were used as benchmarks to assess the performance of the full-waveform approximations.

 \subsection{Regularized Source \( \mathcal{S} \) Confined to a Single Point}
\label{sec:num-s}

When the volumetric source \( s \), or its regularized version \( \mathcal{S} \), is confined to a single point, the spatial integral in the primary Green's formula~\eqref{eq:gr-prime} and its equivalent frequency-domain expression is eliminated. The resulting formula thus describes the action of the causal Green's function on a nonphysical source confined to a single spatial point. This section evaluates the accuracy of the full-waveform approximation in reproducing this analytical action, focusing on comparisons with the exact solution obtained from the primary Green's formula~\eqref{eq:gr-prime}.\footnote{For this simplified and nonphysical numerical experiment using a $d$-dimensional radiation source $s$ confined to a single point, the quantities will have units of $\text{m}^{-\text{d}}$ multiplied by the units of physical quantities associated with a finite-sized source. Integrating over a set of sampled points distributed across a finite volume of the source then yields a wavefield with physically meaningful units.}

The regularized source \( \mathcal{S} \) is a band-limited version of the radiation source \( s \), adapted to the regular space-time grid used for implementing Algorithm~\ref{alg:1}, and supports frequencies up to a prescribed cutoff frequency \( f_{\text{max}} \). We expect the full-waveform approximation of \( \mathcal{S} \) to match the analytical action of the Green's function on \( s \) at all frequencies below \( f_{\text{max}} \). In this section, we examine this agreement in the case where the source is confined to a single spatial point.

For a regularized source confined to an arbitrary single point \( \bx_0 \), the regularized mass source satisfies Eq.~\eqref{eq:mass-source-regularized}, and is discretized as
\begin{align}
\label{eq:mass-s-dis-point}
\bar{\mathcal{S}}_m^{\text{point}}(\bX, \bt) \approx \Delta t \, \bar{\delta}_b(\bX - \bx_0) \sum_{\bt' = 1}^{\bt} \mathcal{S}(\bx_0, \bt'),
\end{align}
which represents a nonphysical quantity due to the omission of the spatial integration.

\subsubsection{Experiment}

\begin{figure} 
   \centering
	\subfigure{\includegraphics[width=0.35\textwidth]{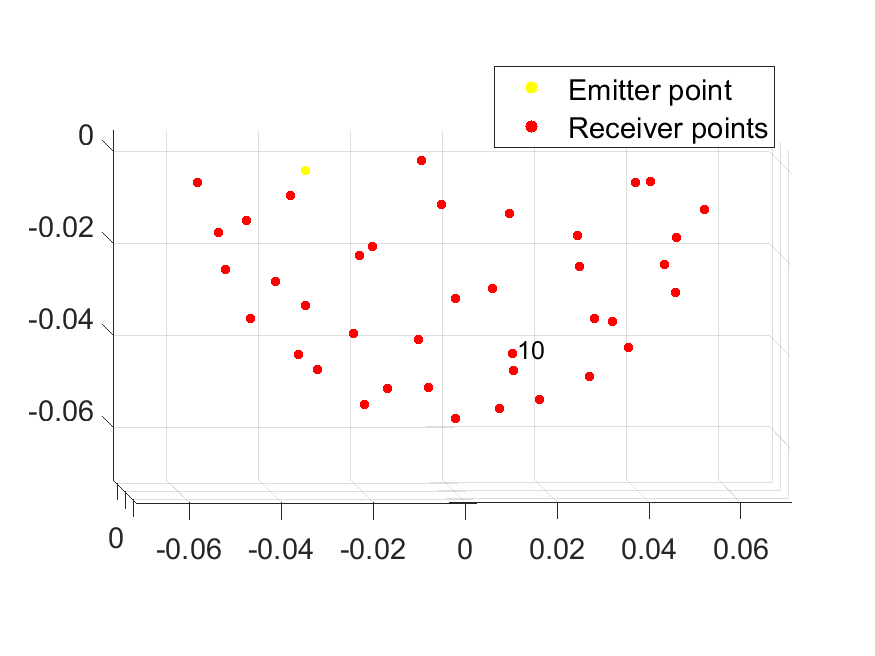}
   \label{fig:1a} }\\
	\subfigure{\includegraphics[width=0.35\textwidth]{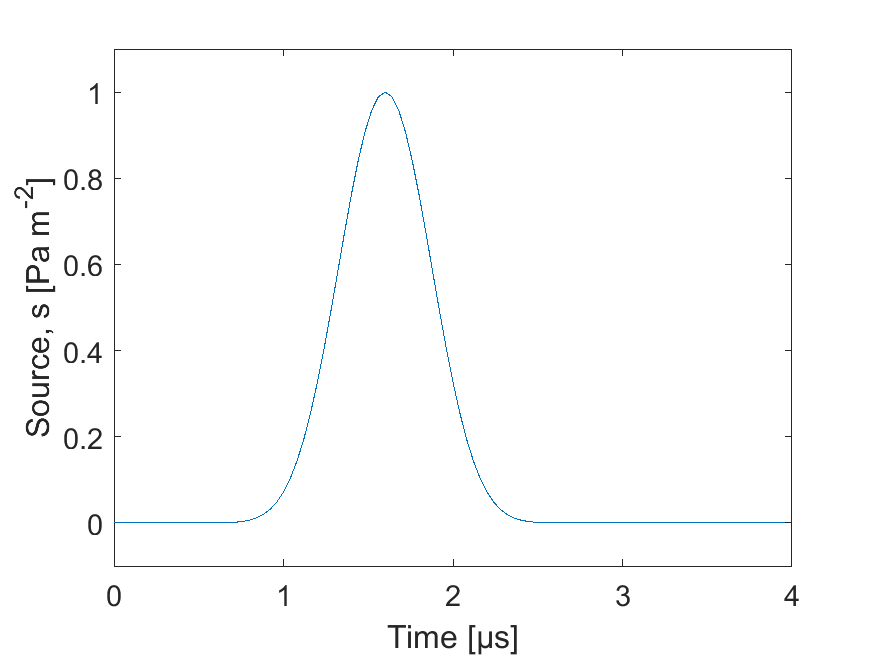}
   \label{fig:1b} }
   	\subfigure{\includegraphics[width=0.35\textwidth]{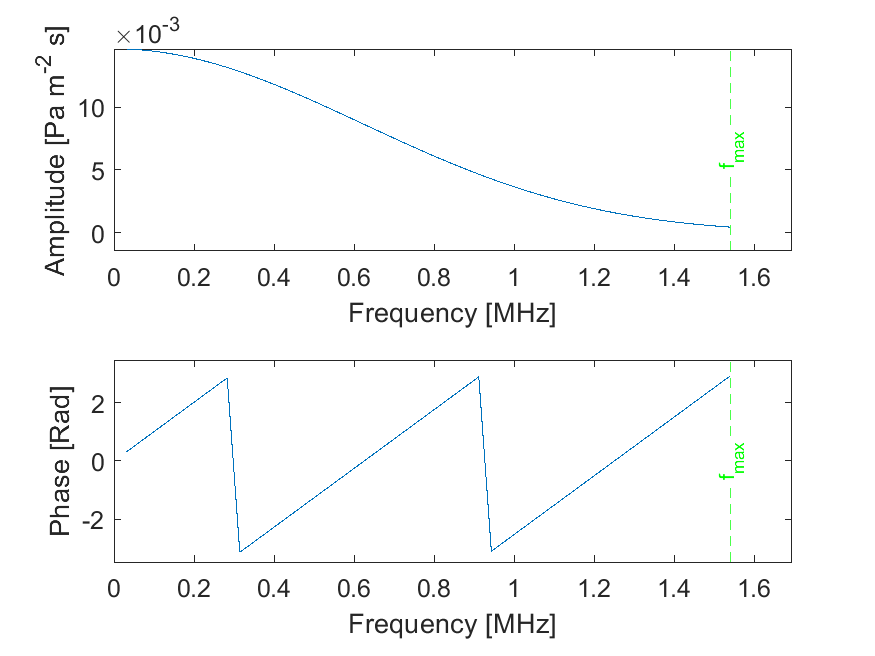}
   \label{fig:1c} }
	\caption{(a) A single emitter point and 39 receiver points arranged on a hemisphere with a radius of 5.6~cm, (b) the source pulse, $s$, represented in the time domain, (c) the source pulse, $s$, shown in the frequency domain, decomposed into amplitude and phase components.}
\end{figure}

A set of 40 points were sampled on a hemisphere centered at the origin of the Cartesian coordinate system with a radius of 5.6~cm. A point, located at $\bx = \big[0, -5.6\big]$~cm, was used as an emitter, while the pressure wavefield was recorded in time on the remaining receiver points. Figure \ref{fig:1a} illustrates the emitter and receiver points, marked in yellow and red, respectively. The source pulse in the time domain is shown in Figure \ref{fig:1b} and is quantified in terms of $s$.

For the point source investigated in this section, a more detailed analysis was performed by comparing results in the frequency domain. Algorithm \ref{alg:1}, with a mass source discretized using Eq.~\eqref{eq:mass-s-dis-point}, was employed to perform a time-domain full-waveform approximation of the action of the causal Green's function on the point source. The time-domain source pulse is depicted in Figure \ref{fig:1b}, while its spatial distribution is represented in Figure \ref{fig:1a}. The approximated, nonphysical wavefield was recorded in time at all receiver positions and subsequently transformed into the frequency domain. This transformation was performed at 50 equidistant discretized frequencies within the range $\big[1/50, 1 \big] \times f_\text{max}$.

For the analytical approach, the source pulse shown in Figure \ref{fig:1b} was transformed into the frequency domain. Its frequency-domain representation, decomposed into amplitude and phase, is shown in Figure \ref{fig:1c}. The action of the frequency-domain Green's function on the frequency-domain source pulse was then calculated at the selected discretized frequencies. In Figure \ref{fig:1c}, the green vertical line indicates $f_\text{max}$, the maximum frequency supported by the computational grid for the full-waveform approximation.

\subsubsection{Results}

\begin{figure} 
\centering
\subfigure[]{\includegraphics[width=0.40\textwidth]{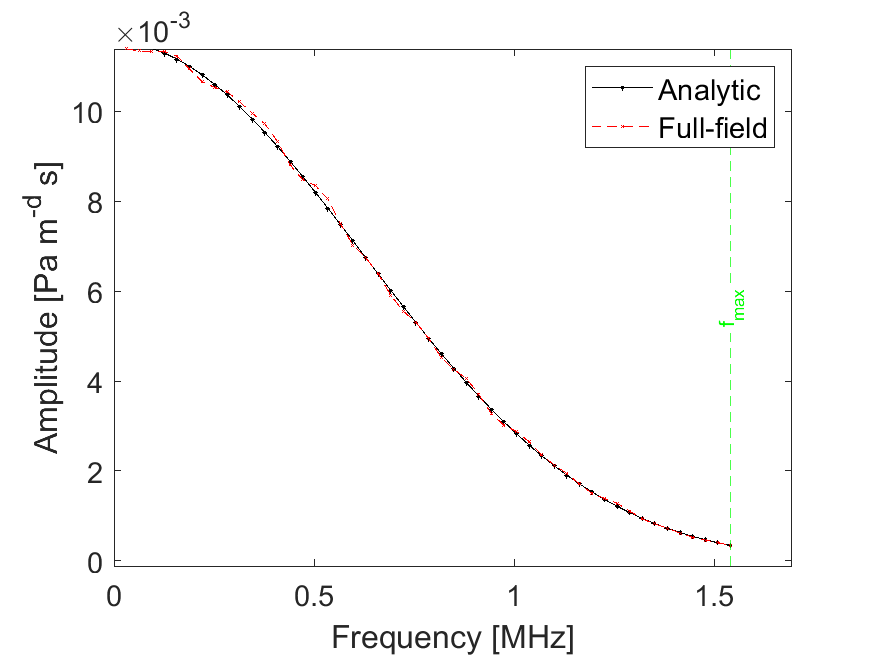}
\label{fig:2a}  }
\subfigure[]{\includegraphics[width=0.40\textwidth]{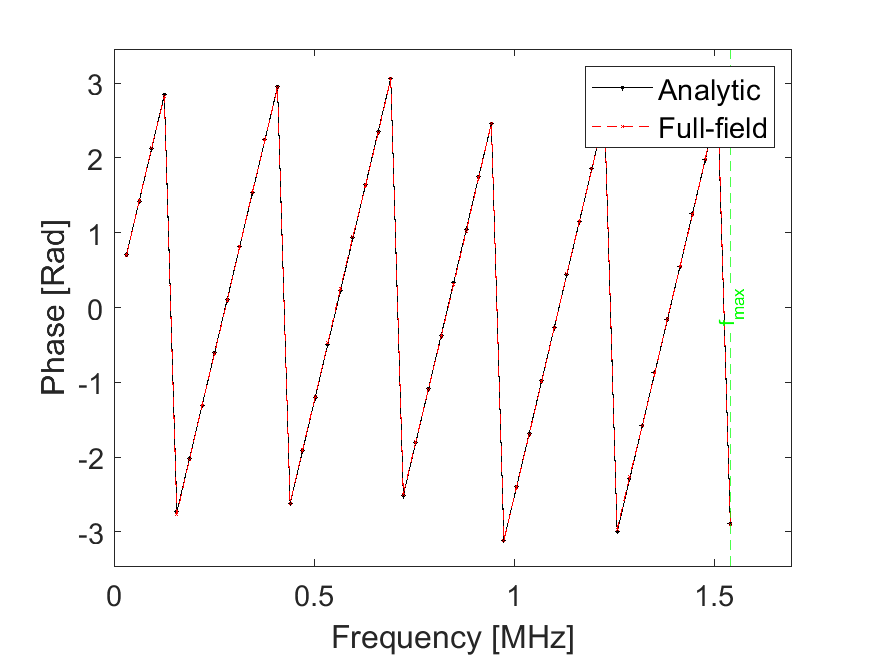}
\label{fig:2b} }
\subfigure[]{\includegraphics[width=0.40\textwidth]{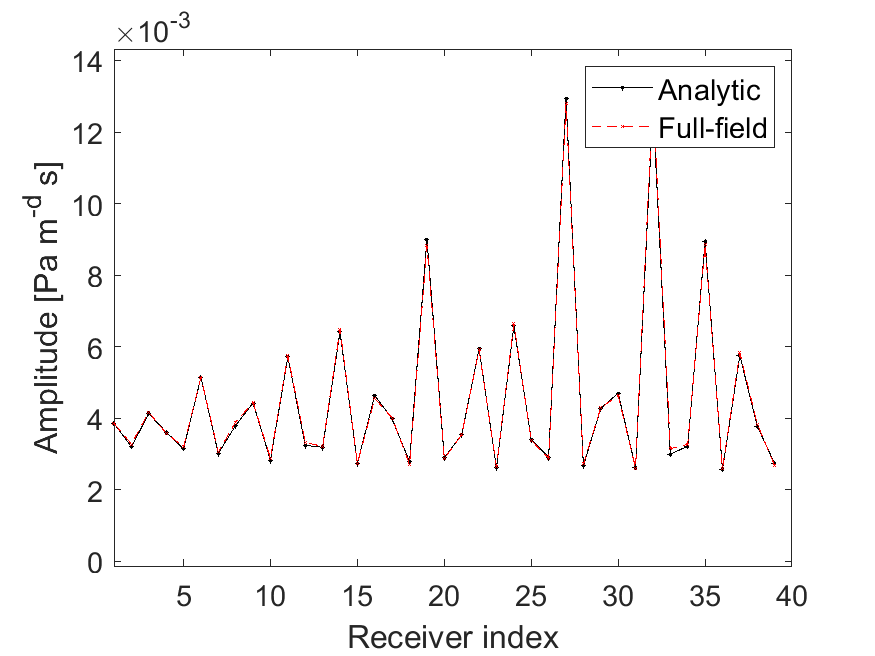}
\label{fig:2c}  }
\subfigure[]{\includegraphics[width=0.40\textwidth]{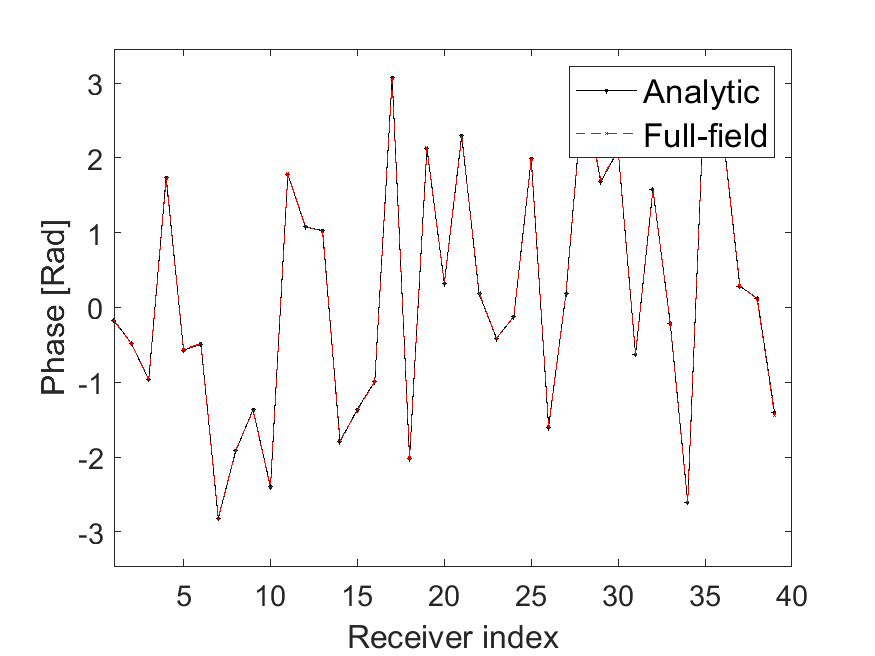}
\label{fig:2d} }
\caption{(a) Amplitude and (b) phase of the action of the causal Green's function on a source defined at a single point, evaluated at receiver 10 across all chosen frequencies. (c) Amplitude and (d) phase of the action approximated at a single frequency of 1~MHz, evaluated across all receiver points.}
\end{figure}

Figure \ref{fig:2a} presents the amplitudes recorded at all selected sampled frequencies on Receiver 10 (as shown in Figure \ref{fig:1a}). The amplitudes calculated analytically using the causal Green's function are displayed in black, while those approximated using the full-waveform approach are shown in red. Similarly, Figure \ref{fig:2b} illustrates the phases computed analytically and approximated via the full-waveform approach at all sampled frequencies on Receiver 10. These figures demonstrate a strong agreement between the analytical solution and the full-waveform approximation for representing the action of the causal Green's function on a point source.

Furthermore, Figures~\ref{fig:2c} and~\ref{fig:2d} depict the amplitude and phase, respectively, of the wavefield at a single frequency of 1~MHz, evaluated across all receiver positions. These results further confirm the agreement between the analytical solution using the radiation source \(s\), confined to a single point, and the full-waveform approximation of its regularized counterpart \(\mathcal{S}\).

\subsection{Monopole Source \( u^{\textbf{n}} \) Supported on a Disk Surface}
\label{sec:re-mo}

This section evaluates the accuracy of a full-waveform approximation of a monopole source, \( \rho_0 \frac{\partial u^{\bn}}{\partial t} \), supported on a surface. The analytical benchmark is provided by the monopole integral formula~\eqref{eq:gr-monopole}, computed using the open-source \emph{Field II} toolbox~\cite{Jensen1,Jensen2}. The full-waveform simulation is performed using Algorithm~\ref{alg:1}, with the mass source discretized according to Eq.~\eqref{eq:mass-u-dis}.

\subsubsection{Experiment}

As discussed in Section~\ref{sec:sum}, the monopole integral formula~\eqref{eq:gr-monopole} is derived under the rigid-baffle assumption and represents the surface integral of the causal Green's function acting on a monopole source, \( \rho_0 \frac{\partial u^{\bn}}{\partial t} \). Here, \( u^{\bn} \) denotes the normal component of the velocity vector relative to the source surface.

\begin{figure} 
\centering
\subfigure[]{\includegraphics[width=0.40\textwidth]{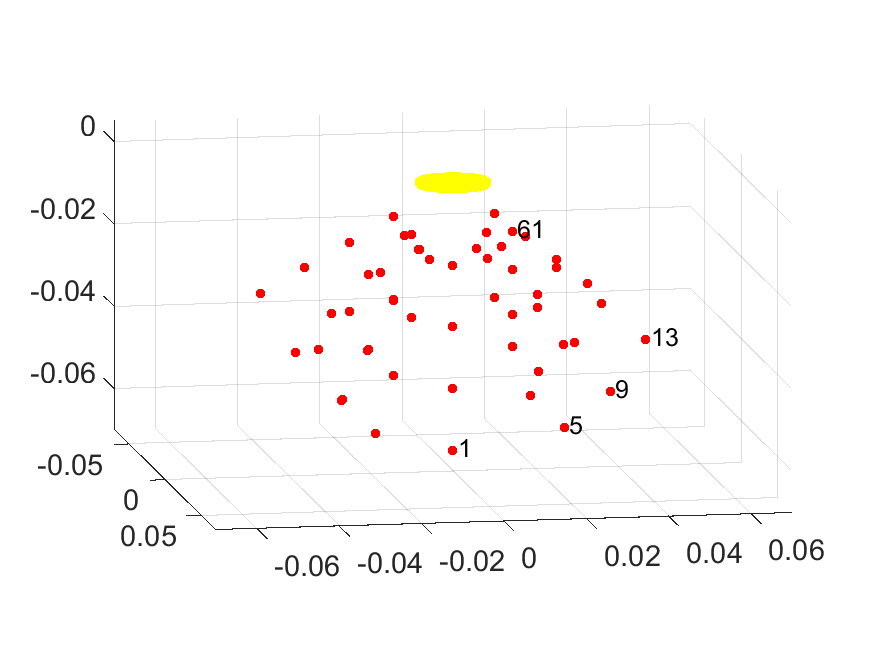}
\label{fig:3a}  }
\subfigure[]{\includegraphics[width=0.40\textwidth]{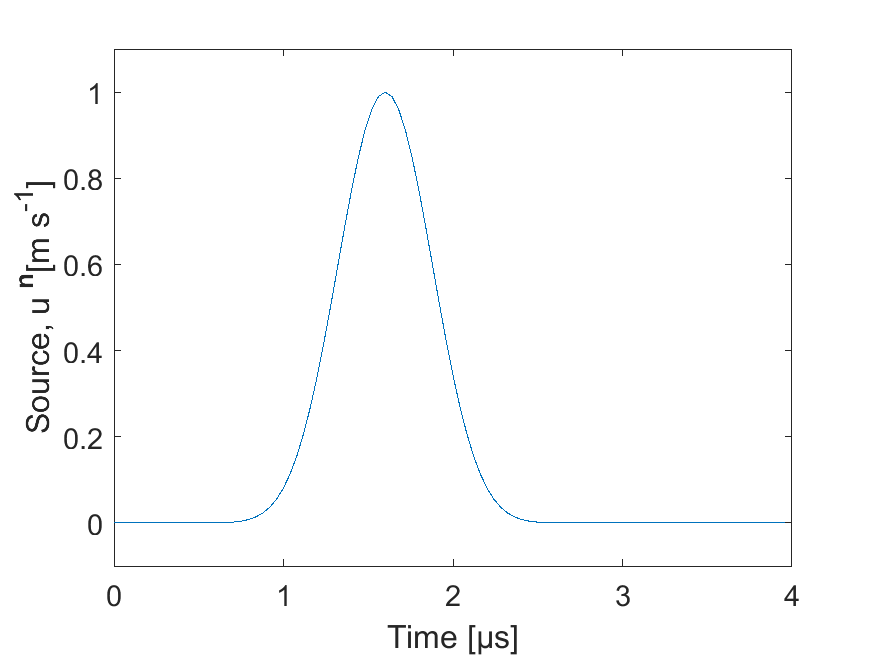}
\label{fig:3b} }
\subfigure[]{\includegraphics[width=0.40\textwidth]{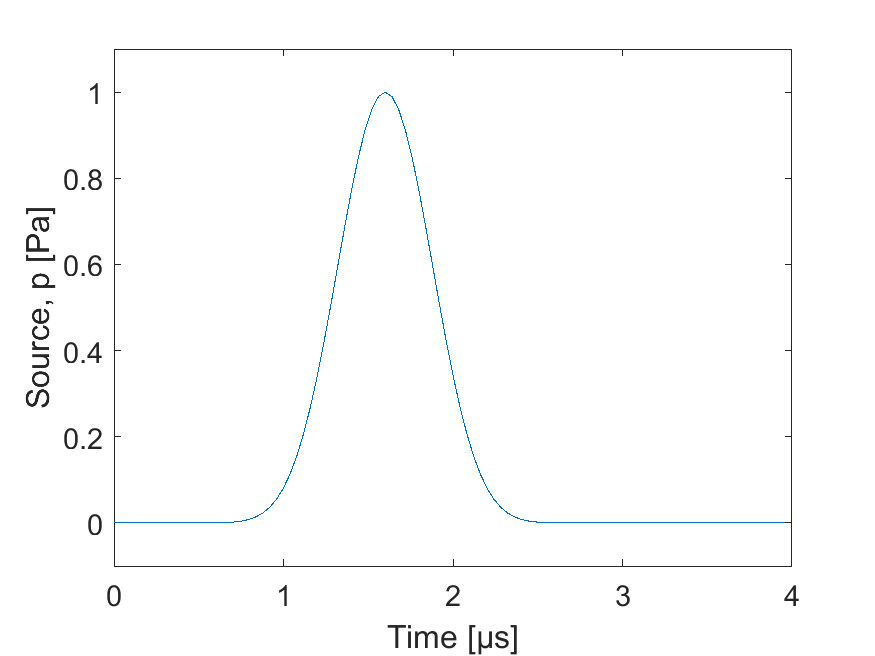}
\label{fig:3c}  }
\caption{(a) A disk-shaped emitter and 64 receiver points, ordered by varying $\theta$, $\varphi$, and $r$. (b) Source pulse $u^{\bn}$ under rigid-baffle conditions. (c) Source pulse $p$ under soft-baffle conditions.}
\end{figure}

In this experiment, the emitter is modeled as a disk-shaped surface with a radius of 8~mm, centered at the origin of the Cartesian coordinate system. The geometry of the emitter surface is illustrated in yellow in Figure~\ref{fig:3a}, while the time-domain source pulse, given in terms of \( u^{\bn} \), is shown in Figure~\ref{fig:3b}.

\begin{enumerate}
   
\item \emph{On-Grid Sampling.}  
In the full-waveform approach, the wavefield is typically approximated at discrete grid points. The grid configuration used in this experiment was described earlier in this section. In contrast, the \textit{Field II} toolbox employs an inherently analytical framework, allowing direct computation of the wavefield at arbitrary spatial positions. For benchmarking purposes, the analytically computed wavefield was sampled at the same grid positions used for the full-waveform approximation.

\item \emph{Off-Grid Sampling.}  
To further evaluate accuracy, a set of 64 off-grid receiver points was used to approximate and record the wavefield. These receiver positions were defined in spherical coordinates \( (r, \theta, \varphi) \), where the radial distances \( r \) from the disk center (origin) were set to \( \{6.5,\, 5,\, 3.5,\, 2\} \) cm. The polar angles \( \varphi \) were chosen as \( \{0,\, \pi/6,\, \pi/4,\, \pi/3\} \), and the azimuthal angles \( \theta \) as \( \{\pi/4,\, 3\pi/4,\, 5\pi/4,\, 7\pi/4\} \). The receiver locations were ordered sequentially by varying \( \theta \), \( \varphi \), and then \( r \).

For example, receiver positions 1--4 correspond to \( r = 6.5\,\text{cm} \) and \( \varphi = 0 \), with the positions distinguished by variations in the azimuthal angle \( \theta \). Similarly, positions 5--8, 9--12, and 13--16 maintain \( r = 6.5\,\text{cm} \) while incrementing \( \varphi \) to \( \pi/6 \), \( \pi/4 \), and \( \pi/3 \), respectively. This pattern is repeated for subsequent sets of receivers, where \( r \) is adjusted to 5, 3.5, and 2~cm for positions 17--32, 33--48, and 49--64, respectively.

For fixed values of \( r \) and \( \varphi \), receiver positions differing only by \( \theta \) (i.e., groups of four consecutive receivers) are symmetric with respect to the emitter disk surface. Consequently, the wavefields at each position within a group are expected to be identical. In particular, for \( \varphi = 0 \) and fixed \( r \), all four values of \( \theta \) correspond to the same spatial point, resulting in 12 redundant receiver indices. These redundant indices are deliberately retained to maintain consistency in the indexing and to avoid confusion in the visualization. The full set of 64 receiver positions, including redundancies, is shown in red in Figure~\ref{fig:3a}, represented in Cartesian coordinates.

\end{enumerate}

\subsubsection{Results}

\begin{figure} 
\centering
\subfigure[]{\includegraphics[width=0.48\textwidth]{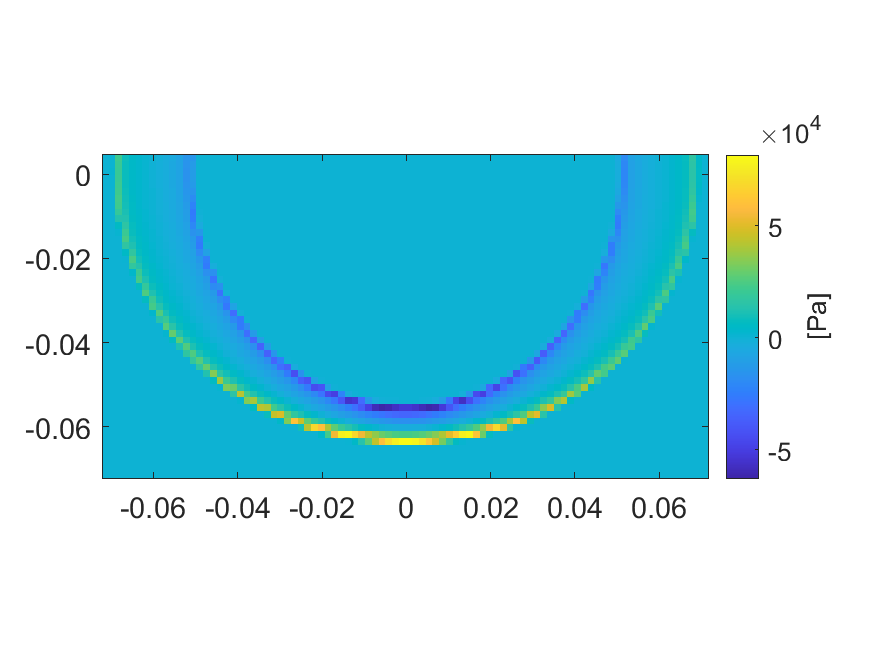}
\label{fig:4a}  }
\subfigure[]{\includegraphics[width=0.48\textwidth]{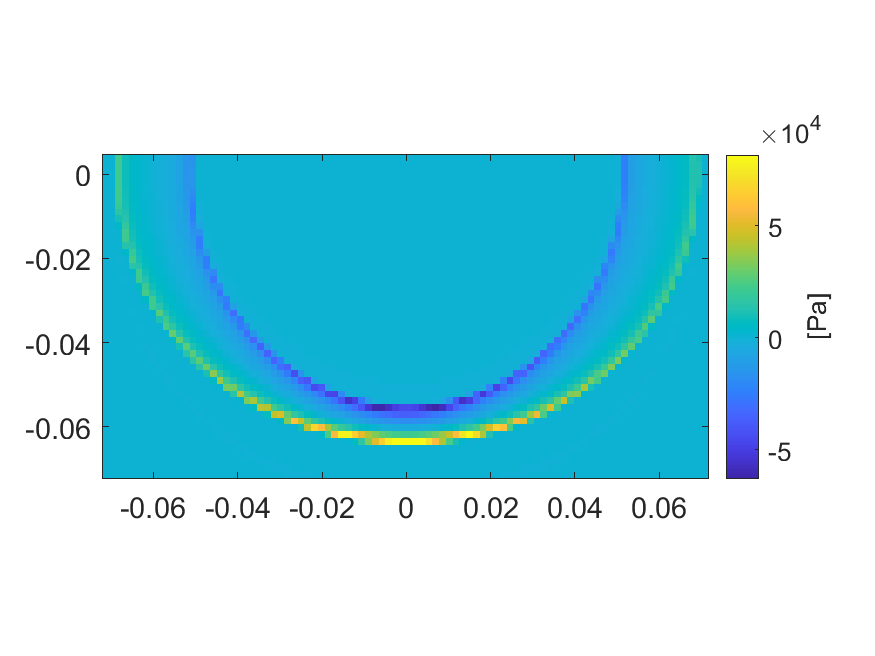}
\label{fig:4b} }
\caption{Wavefields approximated on the plane $\bx^1 = 2.94$~cm at a single time $t = 45\mu\text{s}$, following the excitation of the disk-shaped emitter by the source pulse $u^{\bn}$ (shown in Figure \ref{fig:3b}). The emitter disk's center is positioned at the origin of the Cartesian coordinates, as indicated in yellow in Figure \ref{fig:3a} (not shown here). (a) Analytic solution of Eq.~\eqref{eq:gr-monopole} using \textit{Field II}; (b) Full-waveform approximation using Algorithm~\ref{alg:1} and a mass source discretized via Eq.~\eqref{eq:mass-u-dis}.
}
\end{figure}

Figures~\ref{fig:4a} and \ref{fig:4b} display the wavefields computed analytically using the \textit{Field II} toolbox and approximated using the full-waveform method, respectively. Both wavefields are visualized on grid points located in the plane \( \mathbf{x}^1 = 2.94\,\text{cm} \) at a single time instant, \( t = 45\,\mu\text{s} \). To improve visual clarity, the grid is subsampled by a factor of 4 in both figures. For both the analytical and full-waveform approaches, wavefield recordings began at \( t = 0 \).

The results demonstrate strong agreement between the wavefields obtained analytically and those produced by the full-waveform simulation, validating the accuracy of the latter in approximating the monopole integral formulation.

\begin{figure} 
\centering
\subfigure[]{\includegraphics[width=0.45\textwidth]{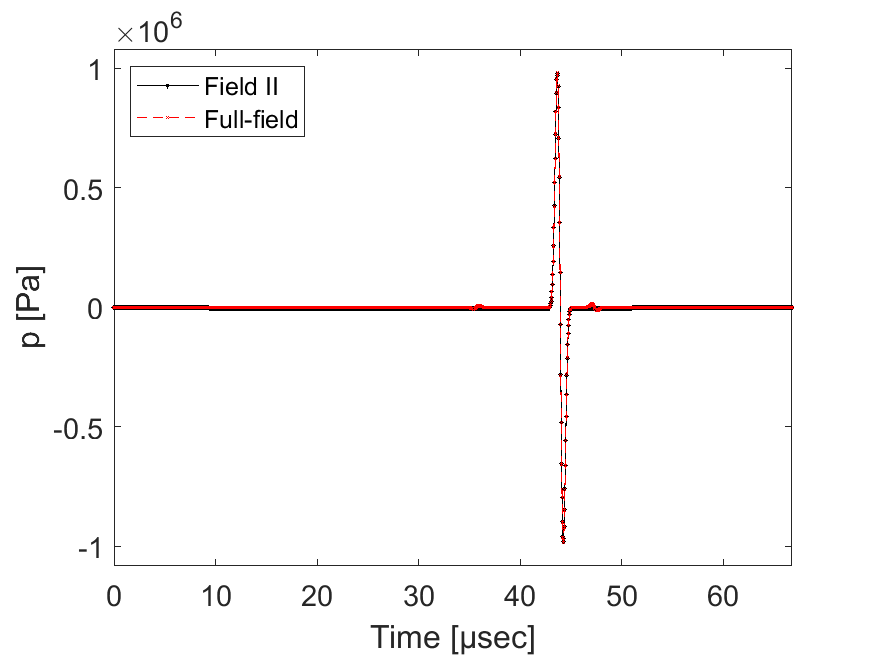}
\label{fig:5a}  }
\subfigure[]{\includegraphics[width=0.45\textwidth]{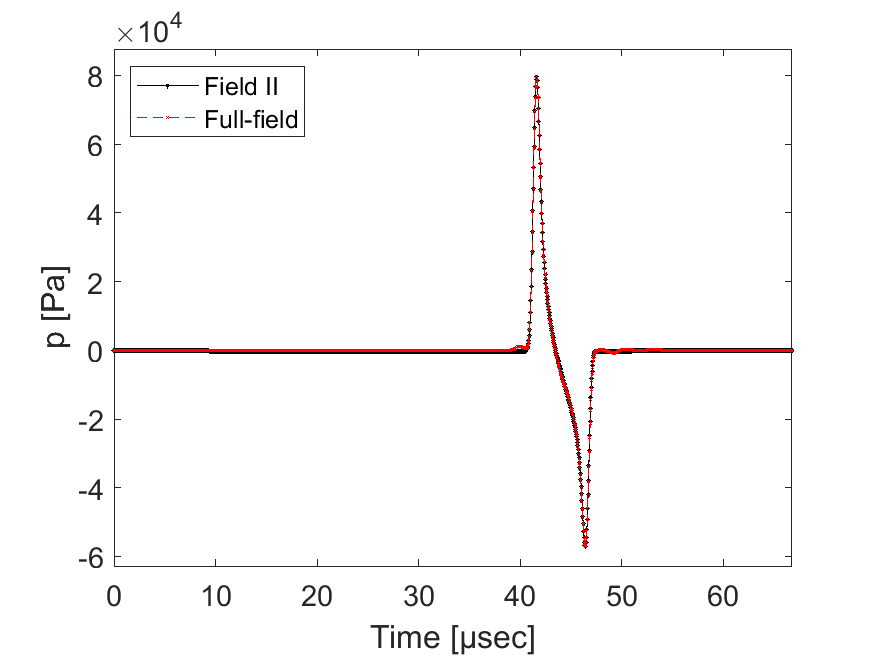}
\label{fig:5b} }
\subfigure[]{\includegraphics[width=0.45\textwidth]{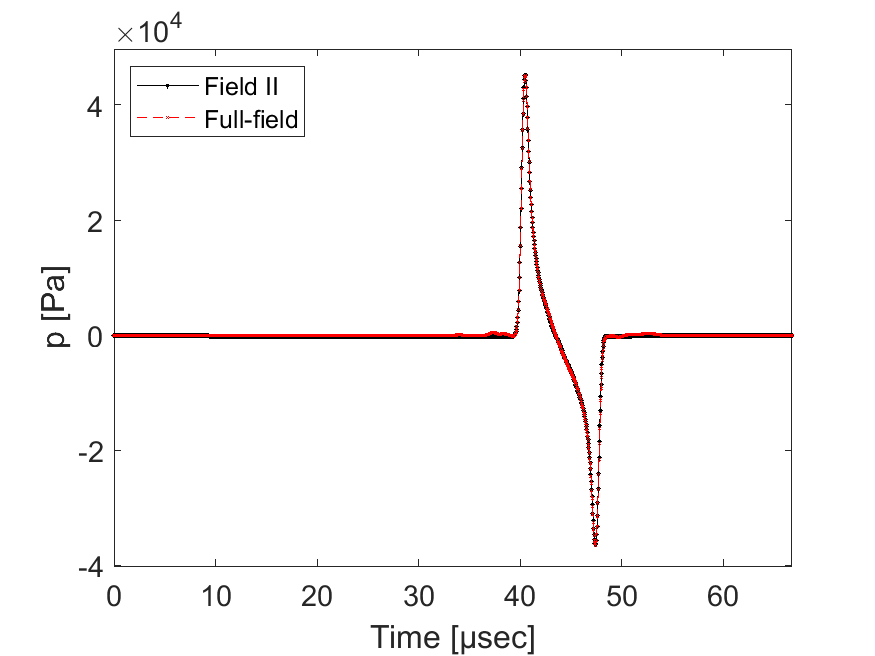}
\label{fig:5c}  }
\subfigure[]{\includegraphics[width=0.45\textwidth]{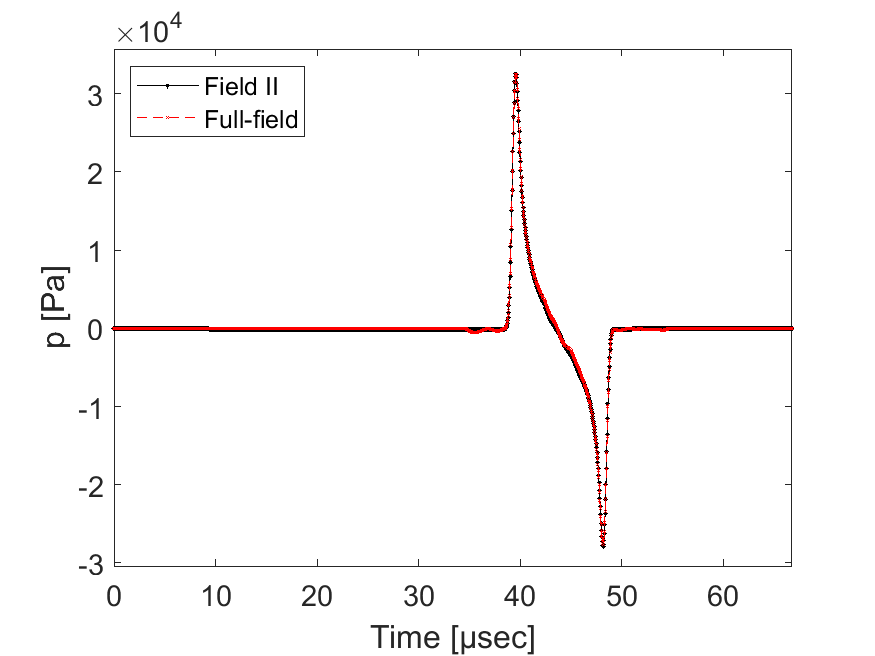}
\label{fig:5d} }
\caption{Wavefields approximated and recorded in time at receiver points following the excitation of the disk-shaped emitter by the source pulse $u^{\bn}$ (shown in Figure \ref{fig:3b}). Receiver points: (a) 1, (b) 5, (c) 9, (d) 13. The monopole integral formula \eqref{eq:gr-monopole} was approximated using both the full-waveform approach and the \textit{Field II} toolbox.}
\end{figure}

Figures~\ref{fig:5a}, \ref{fig:5b}, \ref{fig:5c}, and \ref{fig:5d} illustrate the time-domain wavefields at receiver positions 1, 5, 9, and 13, respectively. These positions, shown in Figure~\ref{fig:3a}, share a common radial distance of \( r = 6.5\,\text{cm} \) and azimuthal angle \( \theta = \pi/4 \), and differ only in the polar angle \( \varphi \). As demonstrated in these figures, the wavefields approximated using the full-waveform method—based on Algorithm~\ref{alg:1} and the discretized mass source formulation of Eq.~\eqref{eq:mass-u-dis}—are in excellent agreement with the corresponding analytical wavefields computed using the \textit{Field II} toolbox.

\begin{figure} 
\includegraphics[width=0.45\textwidth]{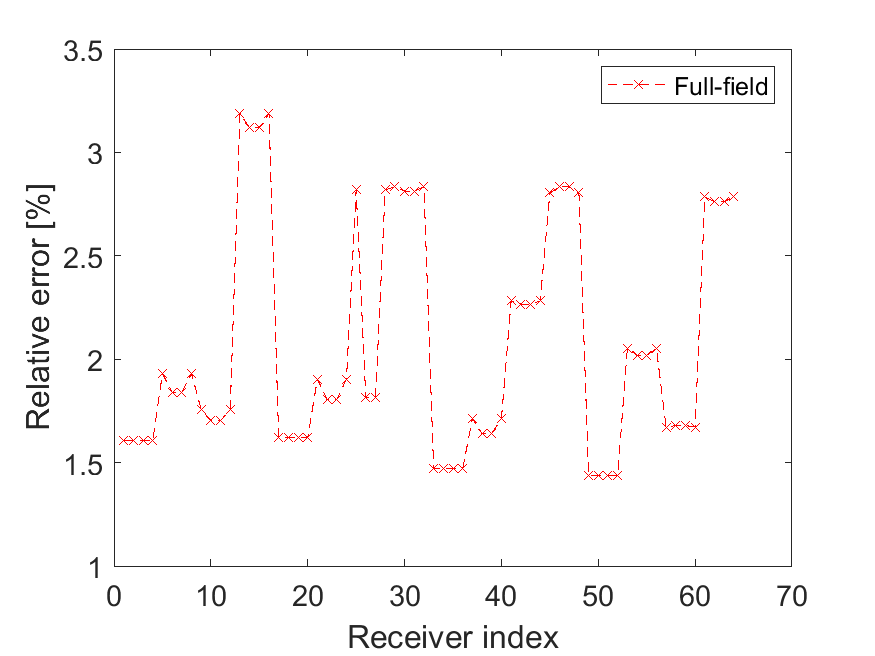}
\caption{Relative error (RE) of the wavefield approximated by the full-waveform approach over time at selected receiver positions, following excitation of the emitter disk by the source pulse \( u^{\bn} \) (shown in Figure~\ref{fig:3b}). The full-waveform results are computed using Algorithm~\ref{alg:1} with a mass source discretized according to Eq.~\eqref{eq:mass-u-dis}. The wavefield obtained from the analytical monopole formula~\eqref{eq:gr-monopole} via the \textit{Field II} toolbox is used as the reference solution. Receiver positions are indexed by varying \(\theta\), \(\varphi\), and \(r\).}
\label{fig:6}
\end{figure}

Figure~\ref{fig:6} presents the relative error (RE) between the full-waveform and analytical wavefields for all 64 receiver positions. Receiver indices are ordered sequentially by varying \( \theta \), \( \varphi \), and \( r \). Recall that within each group of four consecutive receivers (e.g., 1–4, 5–8, ..., 61–64), the parameters \( \varphi \) and \( r \) are held constant, and only \( \theta \) varies. Due to the emitter disk's rotational symmetry, these intra-group positions yield nearly identical wavefield responses, resulting in closely matched RE values. Notably, receiver sets such as 1–4, 17–20, 33–36, and 49–52—corresponding to \( \varphi = 0 \)—are spatially redundant and produce identical RE values. These redundant indices are retained in the plot for completeness and consistency with earlier figures.

This plot demonstrates excellent agreement between the full-waveform approach and the analytical solutions obtained using the \textit{Field II} toolbox, highlighting the accuracy of the wavefield approximations produced by Algorithm~\ref{alg:1} and the discretized mass source defined by Eq.~\eqref{eq:mass-u-dis} in modeling the monopole integral expression~\eqref{eq:gr-monopole}.

\subsection{Dipole Source $ p \textbf{n} $ Supported on a Disk Surface}

This section compares the full-waveform approximation of the dipole integral formula, computed using Algorithm~\ref{alg:1} with a force source discretized according to Eq.~\eqref{eq:force-dis}, to its analytic counterpart given by Eq.~\eqref{eq:gr-dipole1}.

\subsubsection{Experiment}

The emitter disk, shown in yellow in Figure \ref{fig:3a}, is excited by a source pulse \( p \), illustrated in the time domain in Figure \ref{fig:3c}. Under the soft-baffle assumption, the dipole integral formula \eqref{eq:gr-dipole1} is derived, representing the integral of the normal derivatives of the causal Green's function acting on a source \( p \) confined to the surface. For the analytic approximation using the Field\,II toolbox, we reformulate the original dipole integral formula as an integral involving the obliquity-corrected causal Green's functions acting on the surface source, decomposed into far-field and near-field components according to the second line of Eq.~\eqref{eq:gr-dipole1}. The resulting analytic integral was evaluated numerically using the Field\,II toolbox.

The full-waveform approximation of the integral formula \eqref{eq:gr-dipole1} is implemented using Algorithm \ref{alg:1}, with a force source discretized in accordance with Eq. \eqref{eq:force-dis}.

\begin{enumerate}

\item \emph{On-Grid Sampling.}  
As discussed in Section \ref{sec:re-mo}, both the analytic and full-waveform approaches record the wavefield at equispaced sampling points corresponding to the full-waveform approximation grid.

\item \emph{Off-Grid Sampling.}  
Wavefields are approximated and recorded at off-grid receiver positions, shown in red in Figure \ref{fig:3a}, using both the analytic and full-waveform approaches.

\end{enumerate}

\subsubsection{Results}

\begin{figure} 
\centering
\subfigure[]{\includegraphics[width=0.48\textwidth]{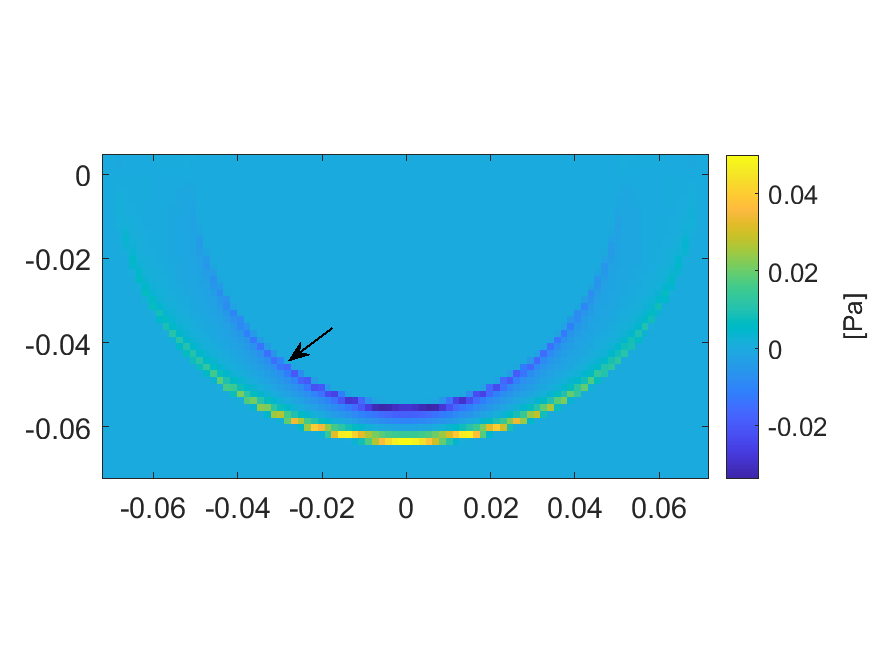}
\label{fig:7a}  }
\subfigure[]{\includegraphics[width=0.48\textwidth]{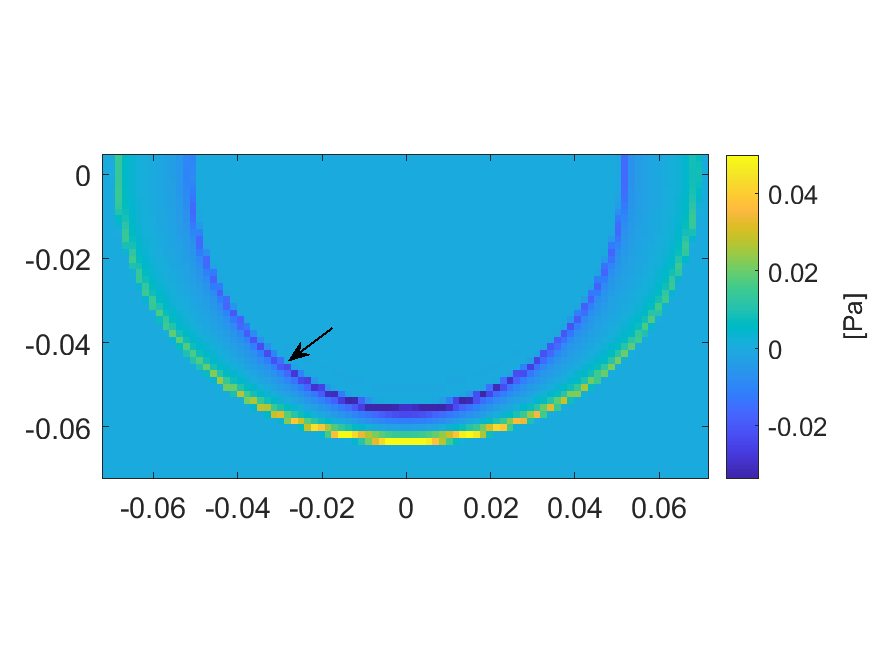}
\label{fig:7b} }
\subfigure[]{\includegraphics[width=0.48\textwidth]{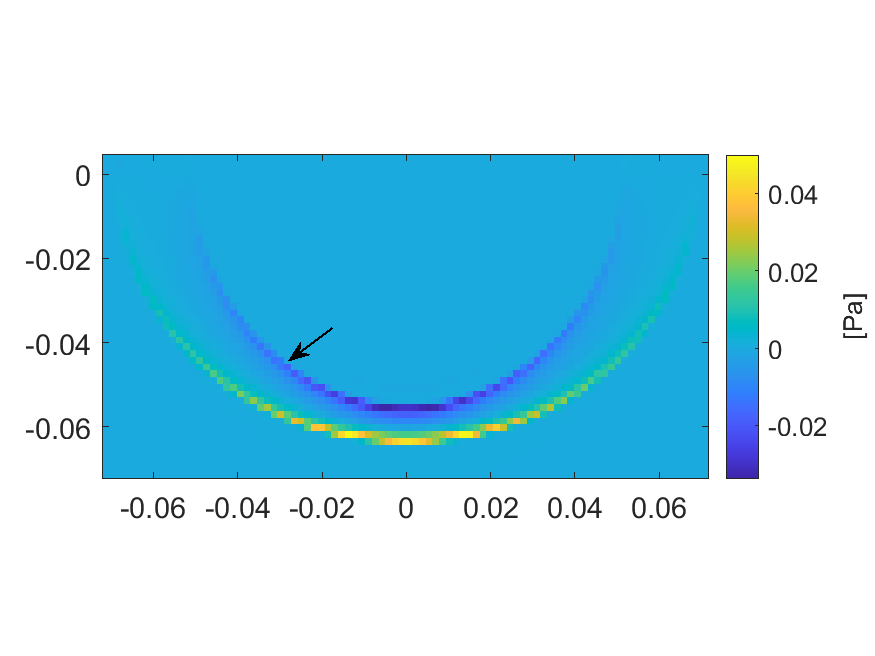}
\label{fig:7c} }
\caption{The wavefield approximated on the plane $\bx^1 = 2.94$~cm and at a time of 45~$\mu\text{s}$ after excitation of the disk-shaped emitter by a source pulse $p$, which is shown in Figure \ref{fig:3c}. The center of the emitter disk is placed at the origin of the Cartesian coordinates, as shown in Figure \ref{fig:3a}. (Not shown here.)  
(a) Analytic solution using \textit{Field II} ,  
(b) Full-waveform approximation using a mass source discretized by Eq. \eqref{eq:mass-p2-dis} added to the equation of continuity,  
(c) Full-waveform approximation using a vector-valued force source defined by Eq. \eqref{eq:force-dis} added to the equation of motion.}
\end{figure}

Figure \ref{fig:7a} shows the wavefields calculated analytically using the dipole integral formula \eqref{eq:gr-dipole1}. Figure \ref{fig:7b} presents the full-waveform approximation of the far-field integral formula \eqref{eq:gr-rsd}, computed using Algorithm \ref{alg:1} with a mass source discretized according to Eq.~\eqref{eq:mass-p2-dis}, which is added to the equation of continuity. As previously discussed, this approximation relies on two key assumptions: (1) the far-field condition \( k x_{\fd} \gg 1 \), and (2) the source is assumed to be omnidirectional, enforced by setting \(\bn \cdot \frac{\bx_{\fd}}{x_{\fd}} = 1\) in Eq.~\eqref{eq:gr-rsd}. In contrast, Figure \ref{fig:7c} illustrates the wavefield approximated using the full-waveform approach in Algorithm~\ref{alg:1}, with a discretized force source defined by Eq.~\eqref{eq:force-dis} added to the equation of motion.

All wavefields are evaluated at grid points lying in the plane \(\bx^1 = 2.94~\text{cm}\) and at a single time instant, \(t = 45~\mu\text{s}\). The grid points are subsampled by a factor of 4 for clarity. For reference, the center of the emitter disk is located at the origin of the Cartesian coordinate system (not shown in the figure).

As seen in Figure~\ref{fig:7b}, the limiting assumptions—particularly the omnidirectionality condition—introduce noticeable discrepancies when compared to the analytic solution in Figure~\ref{fig:7a}. However, Figure~\ref{fig:7c} demonstrates that the full-waveform approximation of the dipole integral formula~\eqref{eq:gr-dipole1}, implemented via the discretized force source term (Eq.~\eqref{eq:force-dis}) incorporated into the equation of motion in Algorithm~\ref{alg:1}, yields a wavefield that closely matches the analytical result.

\begin{figure} 
\centering
\subfigure[]{\includegraphics[width=0.45\textwidth]{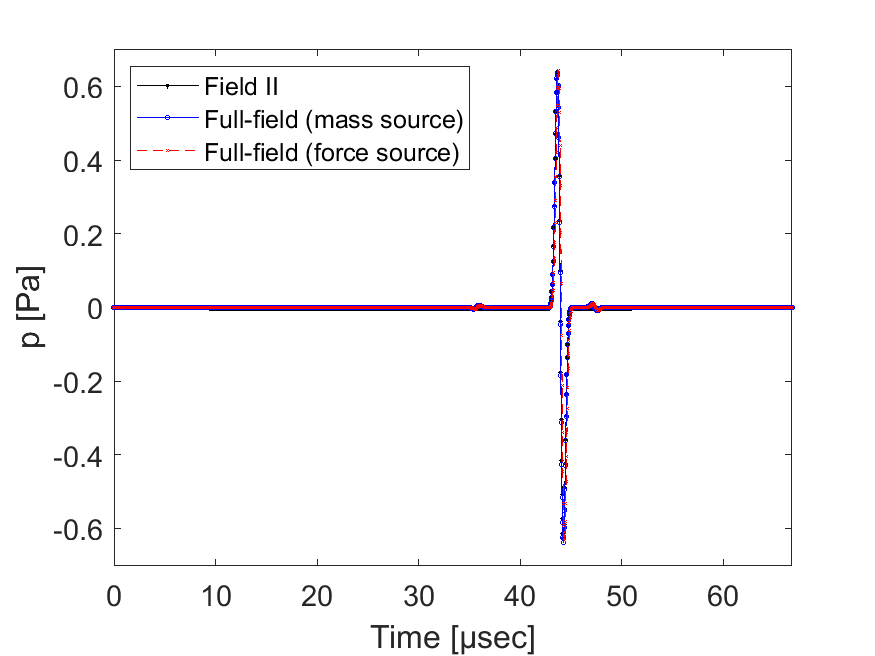}
\label{fig:8a}  }
\subfigure[]{\includegraphics[width=0.45\textwidth]{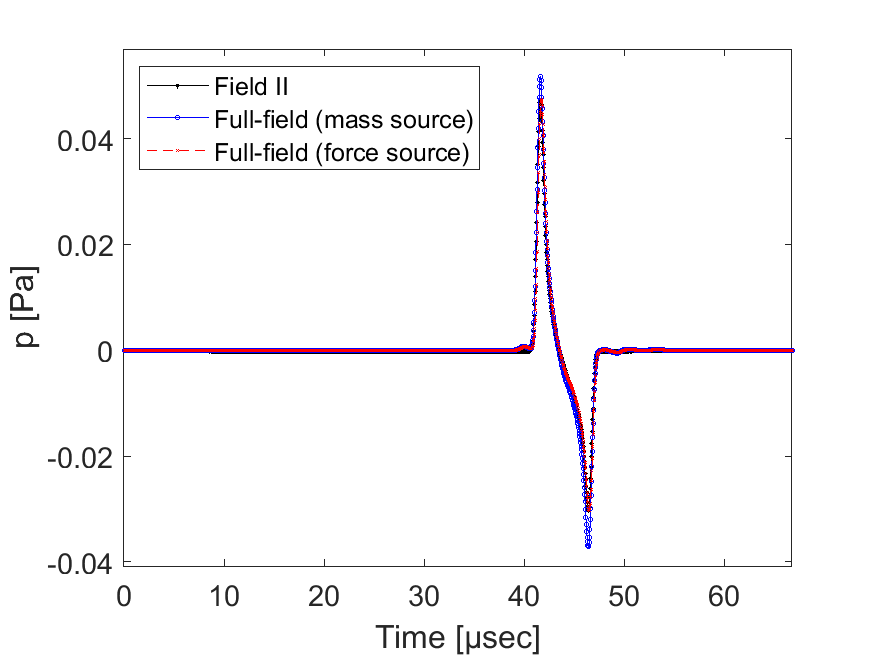}
\label{fig:8b} }
\subfigure[]{\includegraphics[width=0.45\textwidth]{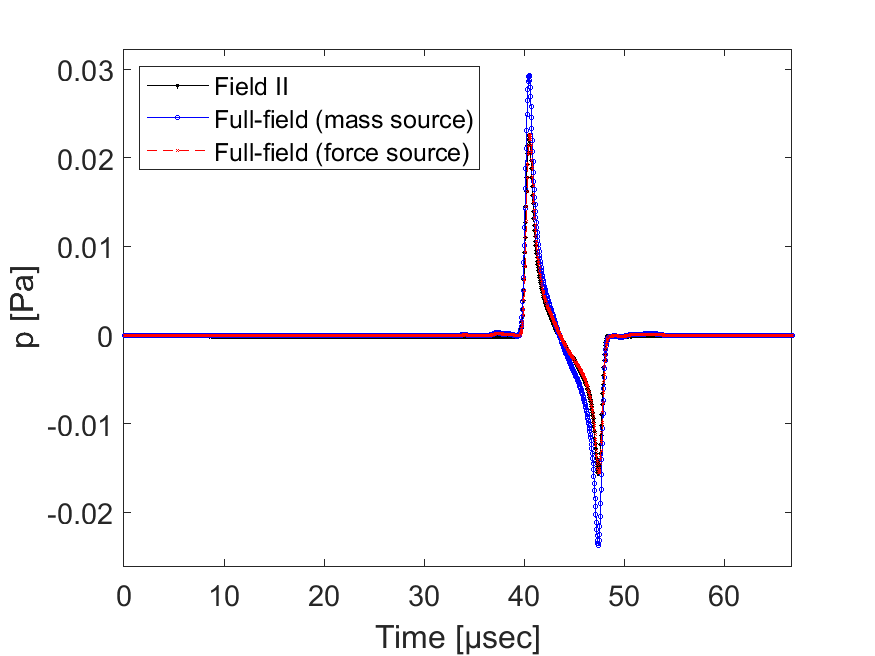}
\label{fig:8c}  }
\subfigure[]{\includegraphics[width=0.45\textwidth]{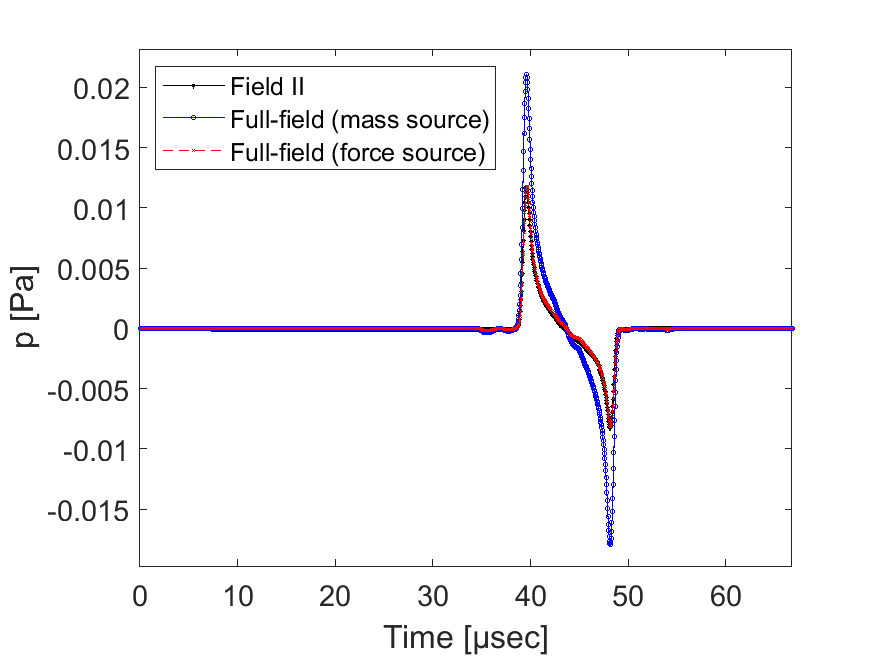}
\label{fig:8d} }
\caption{The wavefield approximated at the receivers over time after the excitation of the disk-shaped emitter by a source pulse $p$, as shown in Figure \ref{fig:3c}. Receiver points: (a) 1, (b) 5, (c) 9, (d) 13. For benchmarking, the \textit{Field II} toolbox is used to analytically calculate the dipole integral formula \eqref{eq:gr-dipole1} (black). The full-waveform approach based on Algorithm \ref{alg:1} is implemented in two ways: (1) using a mass source discretized by Eq. \eqref{eq:mass-p2-dis} added to the equation of continuity (blue), and (2) using a force source discretized via Eq. \eqref{eq:force-dis} added to the equation of motion (red).} \label{fig:8}
\end{figure}

Figures \ref{fig:8a}, \ref{fig:8b}, \ref{fig:8c}, and \ref{fig:8d} show the wavefields approximated and recorded at receiver positions 1, 5, 9, and 13, respectively. These receiver positions are depicted in Figure \ref{fig:3a}. As described in Section \ref{sec:re-mo}, the receiver positions are represented in spherical coordinates for this experiment. Specifically, for the selected receiver positions, the radius \( r \) and the azimuthal angle \( \theta \) are fixed at \( r = 6.5~\text{cm} \) and \( \theta = \pi/4 \), respectively. The positions vary only by the polar angle \( \varphi \), resulting in changes to \( \bn \cdot \bx / x \).

As illustrated in these figures, for all selected receiver positions, the full-waveform approximation of the dipole integral formula~\eqref{eq:gr-dipole1}, computed using Algorithm~\ref{alg:1} with a vector-valued force source discretized via Eq.~\eqref{eq:force-dis}, produces wavefield solutions that closely match those obtained from the analytical formula. However, wavefields approximated using the same algorithm with a mass source discretized using Eq.~\eqref{eq:mass-p2-dis} exhibit large discrepancies compared to the analytical solution calculated using the \textit{Field II} toolbox. These discrepancies are primarily attributed to obliquity effects.

Similarly, Figure \ref{fig:9} presents the wavefields calculated at receiver 61 using the analytical formula \eqref{eq:gr-dipole1} alongside its full-waveform approximations. As shown in figure \ref{fig:3a}, this receiver shares the same $\theta$ and $\phi$ as receiver 13 but is located at \( r = 2 \)~cm. The plot demonstrates strong agreement between the analytical calculation and its full-waveform counterpart obtained using a force source discretized via Eq.~\eqref{eq:force-dis} added to the equation of motion. Moreover, the discrepancy between the analytical calculation and the full-waveform approximation using a mass source discretized via Eq.~\eqref{eq:mass-p2-dis} added to the equation of continuity has increased due to the incorporation of errors arising from the neglect of near-field effects.

\begin{figure} 
\centering
\subfigure{\includegraphics[width=0.45\textwidth]{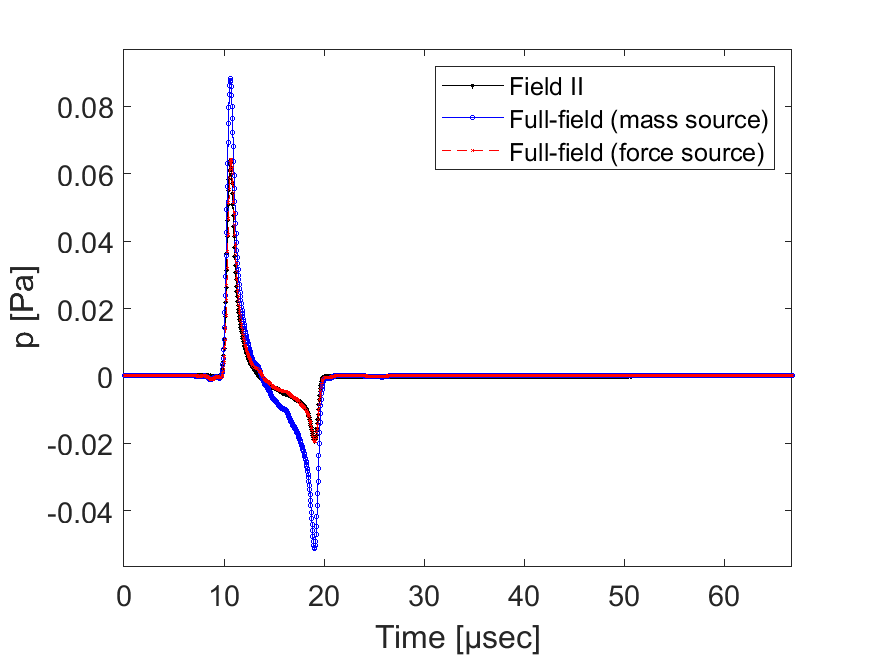}
\label{fig:9}  }
\caption{The wavefield approximated at the receiver point 61 over time after the excitation of the disc-shaped emitter by a source pulse $p$, as shown in Figure \ref{fig:3c}. For benchmarking, the \textit{Field II} toolbox is used to analytically compute the dipole integral formula \eqref{eq:gr-dipole1} (black). The full-waveform approach based on Algorithm \ref{alg:1} is implemented in two ways: (1) using a mass source discretized by Eq. \eqref{eq:mass-p2-dis} added to the equation of continuity (blue), and (2) using a force source discretized by Eq. \eqref{eq:force-dis} added to the equation of motion (red). Compared to receiver point 13 (Fig.~\ref{fig:8d}), the near-field effects have been incorporated in addition to the obliquity effects, leading to large discrepancies between the approximations obtained through the mass and force source definitions.} \label{fig:9}
\end{figure}

\section{Discussion and Conclusion}

Approaches for solving the acoustic wave equation can be broadly categorized into analytic and full-waveform methods. Analytic methods are suitable for homogeneous or weakly heterogeneous (smoothly varying) media, as they can account for refraction effects and singly scattered waves \cite{Jensen1,Jensen2,Javaherian3,Javaherian4,Javaherian-exp}. In contrast, full-waveform approaches are well-suited for handling complex heterogeneities, sharp transitions, and higher-order scattering phenomena in acoustic media \cite{Mast,Tabei,Treeby}.

In certain geophysical applications, acoustic wavelengths are significantly larger than the size of the acoustic aperture, allowing transducers to be approximated as point-like sources or receivers. This assumption, however, does not hold in many biomedical imaging scenarios, where much higher operating frequencies are employed \cite{Xu,Burgholzer1}. Consequently, accurate full-waveform modeling—particularly when amplitude fidelity is required—necessitates accounting for the finite-size effects of acoustic apertures.

This study investigated the equivalence between analytic and full-waveform formulations. In particular, we established how scalar-valued mass sources and vector-valued force sources must be defined, discretized, and incorporated into the full-waveform framework (Algorithm~\ref{alg:1}) to ensure that the resulting solutions are consistent with their corresponding analytical expressions.

To this end, we first demonstrated an equivalence between the analytical action of the causal Green’s function on a volumetric radiation source \(s\) confined to a single point and its full-waveform counterpart obtained using a regularized approximation of this source.

The established equivalence between the analytical and full-waveform approaches for modeling the action of the Green's function on a point source was subsequently extended to model the monopole integral formula~\eqref{eq:gr-monopole}, in which acoustic wavefields are expressed as an integral of the Green's function acting on a surface-supported monopole source, $-\frac{\partial p}{\partial \bn} = \rho_0 \frac{\partial u^{\bn}}{\partial t}$. This source is oriented in the direction normal to the surface.

In addition, it was shown that the dipole integral formula~\eqref{eq:gr-dipole1} is equivalent to an integral of the spatial derivatives of the causal Green's function acting on a dipole source, \( p \bn \), supported on a surface. This integral formula can be reformulated as the second line of Eq.~\eqref{eq:gr-dipole1}, which expresses the wavefield as an integral of obliquity-corrected Green's functions acting on a surface-supported source decomposed into far-field and near-field components.

The far-field formula~\eqref{eq:gr-rsd} is a simplified version of Eq.~\eqref{eq:gr-dipole1}, where the near-field source term, \( \frac{p} { c \, t_{\fd} } \), is neglected. It has been shown that a full-waveform approximation of the far-field dipole formula~\eqref{eq:gr-rsd}, implemented via Algorithm~\ref{alg:1} and using a discretized mass source defined by Eq.~\eqref{eq:mass-p2-dis}, implicitly enforces an additional omnidirectionality assumption, \( \bn \cdot \bx_{\fd} / x_{\fd} = 1 \), thereby treating the dipole source \( p \bn \) as a monopole-like source \( \frac{1}{c} \frac{\partial p}{\partial t} \). This assumption may not hold sufficiently for finite-sized apertures, leading to discrepancies between the approximated wavefield and the analytical solution of the original dipole formula~\eqref{eq:gr-dipole1}.


In contrast, a full-waveform approximation of the dipole integral formula \eqref{eq:gr-dipole1}, implemented using Algorithm \ref{alg:1} with a force source discretized via Eq. \eqref{eq:force-dis}, yielded a wavefield solution that closely matches the analytic solution obtained using the Field II toolbox \cite{Jensen1,Jensen2}. This agreement is further corroborated by the pressure profiles approximated over time, as shown in Figures \ref{fig:8} and \ref{fig:9}.

Leveraging these integral formulations, a reception operator that maps the wavefield in free space onto the receiver surfaces was derived and incorporated into a forward operator within the context of ultrasound tomography. It was demonstrated that the adjoint of this forward operator coincides with a time-reversed variant of the dipole integral formula evaluated on the receiver surfaces.

From the perspective of inverse problems, the observed agreement between the full-waveform approximation of the dipole integral formula~\eqref{eq:gr-dipole1} and its analytic solution is particularly significant. This consistency, along with our numerical results, validates approximating the adjoint and time-reversal operators derived in Section~\ref{sec:adjoint}—a time-reversed variant of the dipole integral formula—via a full-waveform approach. Both the adjoint and time-reversal formulations utilize a time-reversed variant of Algorithm~\ref{alg:1} combined with a force source defined by a time-reversed variant of Eq.~\eqref{eq:force-dis}.

The adjoint or time-reversal operators act on pressure data measured over time on a boundary surface, or on the residual function when employed within iterative frameworks, such as error minimization algorithms~\cite{Dean-Ben,Datchev,Arridge,Javaherian1} or \textit{Neumann series} iterations~\cite{Qian,Haltmeier2}. The next step involves a comprehensive evaluation of the derived full-waveform approximation, explicitly accounting for the finite spatial extent of receivers in practical measurement configurations. The corresponding forward-adjoint operator, incorporating the finite sizes and the analytical angular sensitivity of the receivers, will be compared against established time-reversal and adjoint operators reported in the literature. Particular emphasis will be placed on assessing its integration within adjoint-assisted optimization techniques and the \textit{Neumann-series} framework for solving inverse problems arising in biomedical acoustics, especially those relying on accurate amplitude approximations.

\end{document}